\documentstyle[12pt]{article}

\input amssym.def
\input amssym.tex

\newcommand{\eop}{\bigstar}  

\newcommand{\card}[1]{{\vert #1 \vert} }
\newcommand{\norm}[1]{{\card{\card{#1}}}}

\newcommand{\forces}{\Vdash}

\newcommand{\initial}{\vartriangleleft}
\newcommand{\initialeq}{\vartrianglelefteq}

\newcommand{\ap}{{\rm ap}}

\newcommand{\Dom}{{\rm Dom}}
\newcommand{\Ev}{{\rm Ev}}
\newcommand{\Rang}{{\rm Rang}}

\newcommand{\cf}{{\rm cf}}

\newcommand{\llg}{{\rm lg}}
\newcommand{\md}{{\rm md}}

\newcommand{\implies}{\Longrightarrow}

\newenvironment{Proof}{\noindent{\bf Proof.}}{\par\bigskip} 

\newtheorem{THEOREM}{Theorem}[section]

\newtheorem{Conclusion}[THEOREM]{Conclusion}

\newtheorem{LEMMA}[THEOREM]{Lemma}

\newtheorem{Sublemma}[THEOREM]{Sublemma}

\newtheorem{Main Theorem}[THEOREM]{Main Theorem}
\newenvironment{main Theorem}{\begin{Main Theorem}} 
{\end{Main Theorem}}

\newtheorem{Theorem}[THEOREM]{Theorem}

\newtheorem{Definition}[THEOREM]{Definition}

\newtheorem{Conventions}[THEOREM]{Conventions}

\newtheorem{Main Definition}[THEOREM]{Main Definition}
\newenvironment{main definition}{\begin{Main Definition}}
{\end{Main Definition}}

\newtheorem{Lemma}[THEOREM]{Lemma}

\newtheorem{Notation}[THEOREM]{Notation}

\newtheorem{Convention}[THEOREM]{Convention}

\newtheorem{Note}[THEOREM]{Note}

\newtheorem{Observation}[THEOREM]{Observation}

\newtheorem{Remark}[THEOREM]{Remark}

\newtheorem{Main Fact}[THEOREM]{Main Fact}
\newenvironment{main Fact}{\begin{Main Fact}}{\end{Main Fact}}

\newtheorem{Fact}[THEOREM]{Fact}

\newtheorem{Subfact}[THEOREM]{Subfact}

\newtheorem{Claim}[THEOREM]{Claim}

\newtheorem{Main Claim}[THEOREM]{Main Claim}
\newenvironment{main claim}{\begin{Main Claim}}{\end{Main Claim}}

\newtheorem{Corrolary}[THEOREM]{Corrolary}

\newtheorem{Subclaim}[THEOREM]{Subclaim}

\newtheorem{Corollary}[THEOREM]{Corollary}

\newtheorem{Proposition}[THEOREM]{Proposition}

\newtheorem{Discussion}[THEOREM]{Discussion}

\newenvironment{Proof of the Subfact}
{\noindent{\bf Proof of the Subfact.}}{\par\bigskip}

\newenvironment{Proof of the Theorem}
{\noindent{\bf Proof of the Theorem.}}{\par\bigskip}

\newenvironment{Proof of the Conclusion}
{\noindent{\bf Proof of the Conclusion.}}{\par\bigskip}

\newenvironment{Proof of the Observation}
{\noindent{\bf Proof of the Observation.}}{\par\bigskip}

\newenvironment{Proof of the Fact}
{\noindent{\bf Proof of the Fact.}}{\par\bigskip}

\newenvironment{Proof of the Lemma}
{\noindent{\bf Proof of the Lemma.}}{\par\bigskip}

\newenvironment{Proof of the Sublemma}
{\noindent{\bf Proof of the Sublemma.}}{\par\bigskip}

\newenvironment{Proof of the Claim}
{\noindent{\bf Proof of the Claim.}}{\par\bigskip}

\newenvironment{Proof of the Subclaim}
{\noindent{\bf Proof of the Subclaim.}}{\par\medskip}

\newenvironment{Proof of the Main Claim}
{\noindent{\bf Proof of the Main Claim.}}{\par\bigskip}

\catcode`\@=11
\def\@begintheorem#1#2{\rm \trivlist \item[\hskip \labelsep{\bf #1\ #2.}]}
\def\@opargbegintheorem#1#2#3{\rm \trivlist
      \item[\hskip \labelsep{\bf #1\ #2\ (#3).}]}
\catcode`\@=12

\newcommand{\elementary}{\prec}

\newcommand{\Bbf}{\Bbb}
\newcommand{\vartrianglelefteq}{\trianglelefteq}

\newcommand{\Max}{{\rm Max}}

\newcommand{\into}{\rightarrow}

\newcommand{\rest}{\upharpoonright}  

\newcommand{\deq}{\buildrel{\rm def}\over =}

\newcommand{\CC}{{\cal C}}
\newcommand{\DD}{{\cal D}}

\newcommand{\HH}{{\cal H}}
\newcommand{\II}{{\cal I}}
\newcommand{\JJ}{{\cal J}}
\newcommand{\KK}{{\cal K}}

\newcommand{\TT}{{\cal T}}

\newcount\skewfactor
\def\mathunderaccent#1#2 {\let\theaccent#1\skewfactor#2
\mathpalette\putaccentunder}
\def\putaccentunder#1#2{\oalign{$#1#2$\crcr\hidewidth
\vbox to.2ex{\hbox{$#1\skew\skewfactor\theaccent{}$}\vss}\hidewidth}}
\def\name{\mathunderaccent\tilde-3 }

\begin{document}

\title{On the existence of universal models}

\author{Mirna D\v zamonja\\
School of Mathematics\\
University of East Anglia\\
Norwich, NR4 7TJ, UK\\
\scriptsize{M.Dzamonja@uea.ac.uk}\\
\scriptsize{http://www.mth.uea.ac.uk/people/md.html}\and
Saharon Shelah\\
Mathematics Department\\
Hebrew University of Jerusalem\\
91904 Givat Ram, Israel\\
\scriptsize{shelah@sunset.huji.ac.il}\\
\scriptsize{http://shelah.logic.at/}}
\date{}

\maketitle
\begin{abstract} Suppose that $\lambda=\lambda^{<\lambda}
\ge\aleph_0$, and we are
considering a theory $T$.
We give a criterion on $T$ which is
sufficient for the consistent existence
of $\lambda^{++}$ universal models of $T$ of size $\lambda^+$ for models
of $T$ of size
$\le\lambda^+$, and is meaningful when $2^{\lambda^+}>\lambda^{++}$.
In fact, we work more generally with abstract elementary
classes.
The criterion 
for the consistent existence of universals
applies to various well known
theories, such as triangle-free graphs and simple theories. 

Having in mind possible applications in analysis, we further observe that for such
$\lambda$, for any fixed $\mu>\lambda^+$ regular with
$\mu=\mu^{\lambda^+}$, it is consistent that $2^\lambda=\mu$
and there is no normed vector space
over ${\Bbf Q}$ of size $<\mu$
which is universal for normed vector spaces
over ${\Bbf Q}$ of dimension $\lambda^+$
under the notion of embedding $h$ which specifies $(a,b)$ such
that $\norm{h(x)}/\norm{x}\in (a,b)$ for all $x$.

{\footnote{In the list of publications of S. Shelah, this is
publication
number 614. Both authors thank
the United States-Israel Binational Science Foundation for a
partial support and various readers of the manuscript.
for their helpful comments.


}
}
\end{abstract}

\baselineskip=16pt
\binoppenalty=10000
\relpenalty=10000
\raggedbottom

\section{Introduction.} We study the existence of universal models
for certain natural theories, which are not necessarily first order.
This paper is self-contained, and it continues Saharon Shelah's
\cite{Sh457} and \cite{Sh500}. 
An example of a theory to which our results can be applied is the
theory of triangle-free graphs,
or any simple theory (in the sense of \cite{Sh93}).
For $T$ a theory with a fixed notion
of an embedding between its models, we say that
a model $M^\ast$ of $T$ is universal for models of $T$ (of size $\lambda$)
if every model $M$
of $T$ of size $\lambda$, embeds into $M^\ast$. We similarly define when a
family
of models is jointly universal for models of size $\lambda$.
More generally, we consider universals in an abstract elementary class, see Definition
\ref{absel}. 

Two well known theorems on
the existence of universal models for first order theories $T$
(see \cite{CK}) are
\begin{enumerate}
\item
Under $GCH$, there is a universal model
of $T$ of cardinality $\lambda$ for every
$\lambda>\card{T}$.

\item
If $2^{<\lambda}=\lambda>\card{T}$, then
there is a universal model
of $T$ of cardinality $\lambda$.
\end{enumerate}
Without the above assumptions, it tends to be hard for a
first order theory to have a universal model, see [Sh 457]
for a discussion and further references.

Although
the problem of the existence of universal models for first-order theories
(i.e. elementary classes of all models of such a theory) 
is the one which has been studied most extensively,
there are of course many natural theories which are not first order.
To approach such questions, we view the problem from
the point of view of abstract elementary classes, which were introduced
in [Sh 88] (in \S\ref{definitions} we recall the definitions),
and in a more specialized form earlier by Bjarni J\'onsson, see \cite{CK}.
Such classes will be throughout denoted by $\KK$, and if $\lambda$ is
a cardinal, the family of elements of $\KK$ which have size $\lambda$
will be denoted by $\KK_\lambda$.

In [Sh 457] S. Shelah introduced the notion of an
approximation family and studied 
abstract elementary classes with a ``simple" 
(here called ``workable", to differentiate
them from simple theories in the sense of \cite{Sh93}) $\lambda$-approximation
family. One of the results mentioned in
[Sh 457] is that for $\lambda$ an uncountable
cardinal satisfying $\lambda=\lambda^{<\lambda}$,
it is consistent that every abstract elementary class $\KK$ which has
a workable $\lambda$-approximation family, has an element of size
$\lambda^{++}$ which is universal for
the elements if $\KK$ which have
size $\lambda^+$, i.e. $\KK_{\lambda^+}$. 
Although the main idea of the proof there was correct, there were many
incorrect details and omissions that made the proof and theorem incorrect as
stated.
In this paper
we give a somewhat different proof of this result, and we also deal
with $\lambda=\aleph_0$. Our results give a precise criterion for a class to be
amenable to the theorem about
consistency of the existence of a small family of models in $\KK_{\lambda^+}$
that are universal
for $\KK_{\lambda^+}$. Among the classes which satisfy this criterion are
the class of triangle-free
graphs under embeddings (as shown in [Sh 457]) or in fact
the elementary class of models of any simple theory,
as shown in \cite{Sh500}.

A complete
definition of a $\lambda$-approximation family $K_{\rm ap}$ is given in
\S\ref{definitions}, but let us try to give an intuitive idea here. The easiest way to look at
this is to say that $K_{\ap}$ is a forcing notion whose generic
gives an element of $\KK_{\lambda^+}$. A natural example is
to take a theory $T$, consider the class
of all its models $N$
of size $\lambda^+$ (with universe a subset
of the ordinals $<\lambda^+$),
and define $K_{\ap}$
as the set of all $M$ of size $<\lambda$ which are an elementary
submodel to some such $N$, the order
being $\elementary$. So, for example, the union of an elementary
chain of elements in  $K_{\ap}$ is an element of $\KK$.

As we wish to use approximation families as forcing notions,
we are led to discuss the closure and the chain condition.
$K_{\ap}$ is said to be $(<\lambda)$-smooth, if
every chain of length $<\lambda$ has a least upper bound.
All $\lambda$-approximation
families considered here satisfy this condition.
There are indications
that such an assumption is necessary for the universality results
we wish to obtain, as if smoothness fails strongly there are
no universals, see [GrSh 174].

As we intend to iterate with ($<\lambda$)-supports, our chain condition
has to be a strong version of $\lambda^+$-cc, so to be
preserved under such iterations. The one we
use is $\ast^\varepsilon_\lambda$ from [Sh 288], which is also
the one used in [Sh 457]. This condition is a weakening
of ``stationary $\lambda^+$-cc". We recall the definition at the
beginning of \S\ref{forcing}.
The question now becomes
which $\lambda$-approximation families
yield such a chain condition. We call such
approximation families workable.
This notion is
defined in \S\ref{definitions}.
In
[Sh 457] it is shown that triangle free graphs and the theory
of an indexed
family of
independent equivalence equations have workable
$\lambda$-approximation families. 

In \cite{Sh500} and elsewhere, S. Shelah expresses the view that
the existence of
universal models has relevance to the general problem of classifying
unstable theories. With this in mind, we can consider a theory as ``simple"
if it has a workable approximation family. In 
\cite{Sh93}, another meaning of ``simplicity" is considered: a theory is
called simple if it does not have the tree property. In \cite{Sh500} it was
shown that complete simple first order theories
of size $<\lambda$ have workable
approximation families in $\lambda$.
This can be understood as showing that
all simple theories behave ``better" in the respect of universality than
the linear orders do, as it is known by \cite{KjSh409} that
when $GCH$ fails, linear
orders can have a universal in only a ``few" cardinals. The hope of
finding dividing lines via the existence of universal models is also realized
for some non-simple theories, as it was shown by S. Shelah in \cite{Sh457}
that some non-simple theories have workable approximation families,
like the triangle-free graphs and the theory of an indexed family of
independent equivalence relations, as simplest prototypes of non-simple
theories. In \cite{Sh500}, S. Shelah
introduced a hierarchy ${\rm NSOP}_n$
for $3\le n\le\omega$ with the intention of encapturing
by a formal notion the class of first order theories
which behave ``nicely" with respect to having universal models.
Our research here continues \cite{Sh457}.

We now give an idea of the proof of
the positive consistency results. Details are
explained in \S\ref{definitions} and
\S\ref{forcing}. The idea is that through a
$(<\lambda)$-supports iteration of $(<\lambda)$-complete forcing
we obtain the situation under which to every workable strong
$\lambda$-approximation family
$K_{\ap}$ there corresponds a tree of elements of $K_{\ap}$.
If $K_{\ap}$ approximates $\KK$ and $\KK$ is nice enough,
then the
models in this tree are organized so that the entire tree can be amalgamated
to a model in $\KK_{\lambda^+}$.
Along the iteration we also make sure that
every element of $\KK_{\lambda^+}$ can be embedded into a model
obtained as the union of one branch of such a tree.
There are $\lambda^{++}$ trees used for every approximation
family, so the universal model obtained has size $\lambda^{++}$.
Every
individual forcing used in the iteration has $\ast^\varepsilon_{\lambda}$,
but the proof of this for $\lambda>\aleph_0$ requires us to introduce an
auxiliary step in the forcing.

In \S\ref{negativex} we give a consistency
result showing that with the same assumptions on $\lambda^+$
as above it is consistent  that
there is no universal normed vector space of size $\lambda^+$, even under a
rather
weak notion of embedding. We note that negative {\em consistency} results relevant to the
universality problem tend to be much easier to obtain than the positive ones,
especially as far as the first order theories are concerned.

We finish this introduction by giving more remarks on related results,
and some conventions used throughout the paper.

The pcf theory of S. Shelah
has proved to be a useful line of approach to the negative aspect of the problem
of universality. This approach
has been extensively applied by
Menachem Kojman and S. Shelah  (e.g. to linear orders \cite{KjSh409}), and later
by each of them separately (M. Kojman on graphs \cite{Kj},
S. Shelah on Abelian groups
\cite{Sh552}
e.g). See 
[Sh 552] for the history and more references. One of the ideas
involved is to use the existence of a club guessing sequence
to prove that no universals exist. A related result
of Mirna D\v zamonja in \cite{Dz} deals with uniform Eberlein compacta, and in 
\cite{Dz2} she shows how the universality axioms presented in this paper can be
applied to that class.
Among the positive universality results,
let us quote
a paper by Rami Grossberg and S. Shelah \cite{GrSh174},
in which it is shown that e.g. the class
of locally finite groups has a universal model in any strong limit of cofinality
$\aleph_0$ above a compact cardinal.
This paper is also the first reference to the consideration of the
universality spectrum as a useful dividing line in model theory.

Further positive consistency results appear e.g in S. Shelah's \cite{Sh 100}
where the consistent existence of a universal linear order at $\aleph_1$
with the negation of $CH$ is shown, and in S. Shelah's \cite{Sh175},
\cite{Sh175A} where the consistency of
the existence of a universal graph at $\lambda$
for which there is $\kappa$
satisfying $\kappa=\kappa^{<\kappa}<\lambda<2^\kappa=\cf(2^\kappa)$,
is proved. The latter result was continued by Alan Mekler in \cite{M},
where \cite{Sh175} was extended to a larger class of models.

Relating to our negative consistency result, the problem of universality
has been extensively studied in functional analysis, most often for classes of
Banach spaces. Probably the earliest result here is one of
Stefan Banach himself in \cite{Ba} in which he
showed that $C[0,1]$ is isometrically universal for separable
Banach spaces.
Another well known result is that of Wies\l aw Szlenk,
showing that there is no universal separable reflexive
Banach space, \cite{Sz}.
Jean Bourgain expanded on these ideas to
build a body of work. The combinatorial approach to the problem of
universality in spaces coming from functional analysis is used in Stevo Todor\v cevi\'c's
\cite{To}.

Model theory as an approach to study of Banach spaces
has been extensively used, for example by Jean-Louis Krivine
in \cite{Kr} and C.~Ward Henson in \cite{He}. See
Jacques Stern's \cite{Stern} for an account on the early history of
this interaction and \cite{Io1}
for a more recent history. Of the work of this area which is being
currently carried on, we mention
a systematic attempt to a classification theory for Banach spaces
by Jos\'e Iovino, see e.g \cite{Io2}, \cite{Io3},
which also give historical
remarks.

\begin{Convention} (1) We make the standard assumption that the family
of forcing names that we use is full, i.e. if
$p\forces``(\exists x)[\varphi(x)]"$, then
there is a name $\name{\tau}$ such that $p\forces``[\varphi(\name{\tau})]"$.

{\noindent (2)} If $\kappa=\cf(\kappa)<\alpha$, we let
\[
S^\alpha_{\kappa}\deq\{\beta<\alpha:\,\cf(\beta)=\kappa\}.
\]

{\noindent (3)} $\chi$ is throughout assumed to be a large enough regular
cardinal. $<^\ast_\chi$ stands for a fixed well ordering of the set of
all sets hereditarily of size $<\chi$, namely
$\HH(\chi)$.

{\noindent (4)} lub stands for the ``least upper bound", i.e. $M$ is the
lub of a set $\cal M$ in the order $\le$ iff it is its unique least upper bound,
which means that $M$ is an upper bound of $\cal M$ and for every $M^\ast$ such that
$(\forall N\in {\cal M})\,[N\le M^\ast]$, we have $M\le M^\ast$.

{\noindent (5)} For a model $M$, we use $\card{M}$ to denote the underlying
set of $M$, and hence $\norm{M}$ to denote the cardinality of $\card{M}$.
\end{Convention}

\section{Approximation families.}\label{definitions}

\begin{Definition}\label{lawful} [Sh 457] Given $\lambda$ an infinite cardinal,
and
$u_1,u_2\subseteq \lambda^+$.

A function $h:\,u_1\to u_2$
is said to be {\em lawful} iff it
is 1-1 and for all $\alpha\in u_1$
we have $h(\alpha)+\lambda=\alpha+\lambda$.
\end{Definition}

\begin{Notation} (1) For $A\subseteq \lambda^+$, let 
\[
\iota(A)\deq\min
\{\delta:\,A\subseteq\delta\,\,\&\,\,\lambda|\delta\}.
\]
If $M$ is a model, we let $\iota(M)\deq\iota(\card{M})$.
 
{\noindent (2)} In the following, we shall use the notation
$M\rest\delta$ for $M\rest\tau_{\iota(M\cap\delta)}\rest\delta$
(the meaning of $\tau$ and $M$ will be described in the following definition).

{\noindent (3)} $\Ev\deq\{2\beta:\,\beta<\lambda^+\}$.

\end{Notation}

\begin{Remark}\label{andreas} The notion of divisibility of ordinals used here is that
$\lambda|\delta$ means that $\delta = \lambda \cdot \xi$ [not $\delta=
\xi \cdot \lambda$] for some $\xi$. The intuition behind the definition of
a lawful function is that one regards
$\lambda^+$ as partitioned into blocks of length $\lambda$, and then a function
is lawful iff it acts by permuting within each block. Then the function $\iota(A)$
simply measures how far the blocks go that meet $A$.
\end{Remark}
 
\begin{Definition}\label{general} [Sh 457] Let $\lambda$ be an infinite
cardinal.
\begin{description}
\item{(1)} Pair $K_{\ap}=(K_{\rm ap},\le_{K_{\ap}})$
is {\em a weak $\lambda$-approximation family} iff
for some (not necessarily strictly)
increasing sequence{\footnote{For the applications mentioned in this paper, in the
following definitions readers can restrict their attention to the
situation of $\tau_i=\tau_0$ for all $i$.}}
\[
\bar{\tau}=\langle\tau_{i}:\,i<\lambda^+\,\,\&\,\,\lambda|i\rangle\
\]
of finitary vocabularies, each of size $\le \lambda$
we have
\begin{description}

\item{(a)} $K_{\ap}$ is a set partially ordered by
$\le_{K_{\ap}}$, and such that
\[
M\in K_{\ap}\implies M\mbox{ is a }\tau_{\iota(M)}\mbox{-model}.
\]

\item{(b)} If $M\in K_{\ap}$, then $\card{M}\in
[\lambda^+]^{<\lambda}$
and $M\le_{K_{\ap}} N\implies M\subseteq N$.

\item{(c)} If $M\in K_{\ap}$ and $\lambda |\delta$, then
$M\!\rest\delta\in K_{\ap}$ and $M\!\rest\delta\le_{K_{\ap}} M$.
Also
{\footnote{The following contradicts the usual notation of model theory of
forbidding empty models, as in such a situation we cannot interpret individual constants.
However, the meaning of $\emptyset$ we use is clear.}},
$\emptyset=M\rest 0 \in K_{\ap}$. If $M,N\in K_{\ap}$ and $\lambda
|\delta$, while $M\le_{K_{\ap}}N$, then $M\rest\delta \le_{K_{\ap}}N\rest\delta$.

\end{description}

\item{(2)} With $K_{\ap}$ as in (1), a function $h$ is said to
be a $K_{\ap}$-{\em isomorphism} from $M$ to $N$
iff $\Dom(h)=M, \Rang(h)=N$
are both in $K_{\ap}$,
and $h$ is a
$\tau_{\iota(M)}$-isomorphism.

\item{(3)} A weak $\lambda$-approximation family
$(K_{\rm ap},\le_{K_{\ap}})$ is said to be a {\em strong
$\lambda$-approximation family} iff in addition to (a)--(c)
above, it satisfies:

\begin{description}
\item{(d)} [Union] Suppose that $i^\ast<\lambda$. 

If $\bar{M}=\langle
M_i:\,i<i^\ast\rangle$
is a $\le_{K_{\ap}}$-increasing sequence in $K_{\ap}$,
then we have that $\bigcup_{i<i^\ast}M_i$ is
an element of $K_{\ap}$, and it is the $\le_{K_{\ap}}$-lub
of $\bar{M}$.

\item{(e)} [End extension/Amalgamation] If $0<\delta<\lambda^+$ is divisible by
$\lambda$, and $M_0, M_1, M_2\in K_{\ap}$ are such that $M_2\rest\delta= M_0
\le_{K_{\ap}} M_1$ and $\card{M_1}\subseteq \delta$, then
$M_1$ and $M_2$ have a $\le_{K_{\ap}}$-upper bound $M_3$ such that
$M_3\rest\delta=M_1$.

If $M_0,M_1,M_2,\delta$ are as above and $M_1, M_2\le M$, then there is
$M_3\le M$ such that $M_3\ge M_1, M_2$ and $M_3\rest\delta=M_1$.

\item{(f)} [Local Cardinality] For $\alpha<\lambda^+$, the set $\{M\in
K_{\ap}:\,\card{M} \subseteq\alpha\}$ has cardinality $\le\lambda$.

\item{(g)} [Uniformity] For $M_1,M_2\in K_{\ap}$,
we call $h:\,M_1\into M_2$ a {\em lawful isomorphism} iff $h$ is a lawful
function and a $K_{\ap}$-isomorphism. We demand
\begin{description}
\item{$(\alpha)$} if $M\in K_{\ap}$ and $h$ is a lawful mapping from
$\card{M}$ onto some $u\subseteq\lambda^+$, then
for some $M'\in K_{\ap}$ we have that $\card{M'}=u$ and $h$
is a lawful $\tau_{\iota\card{M}}$-isomorphism from $M$ onto $M'$.

\item{$(\beta)$} lawful $K_{\ap}$-isomorphisms preserve $\le_{K_{\ap}}$.
\end{description}

\item{(h)} [Density] For every $\beta$ in $\lambda^+$, and
$M\in K_{\ap}$, there is $M'\in K_{\ap}$
such that $M\le_{K_{\ap}} M'$ and $\beta\in \card{M'}$.

\item{(i)} [Amalgamation] Assume $M_l\in K_{\ap}$ for $l<3$ and $M_0
\le_{K_{\ap}} M_l$
for $l=1,2$. Then for some lawful function $f$
and $M\in K_{\ap}$, we have
$M_1\le_{K_{\ap}} M$, the domain of $f$ is $M_2$, the restriction
$f\rest \card{M_0}$
is the identity, and $f$ is a $\le_{K_{\ap}}$-embedding of $M_2$ into $M$,
i.e. $f(M_2)\le_{K_{\ap}} M$. If $M_1\cap M_2=M_0$, we can assume that
$f={\rm id}$.

\end{description}
\end{description}
\end{Definition}

\begin{Remark}\label{concern} (1) There is no contradiction concerning
vocabularies in (g)$(\alpha)$ of
Definition \ref{general}(3): if $K_{\ap}$ is a weak $\lambda$-approximation
family, while $M\in K_{\ap}$ and $h$ is a lawful mapping from
$\card{M}$ onto some $u$, then $\iota(u)=\iota(\card{M})$ (so
saying that $h$ gives rise to a $K_{\ap}$-isomorphism makes sense).

[Why? Letting $\delta\deq \sup(u)$, if $\gamma<\delta$, we can find
$\alpha\in \card{M}$ such that $h(\alpha)\in (\gamma,\delta)$.
Hence
\[
\gamma<\gamma+\lambda\le h(\alpha)+\lambda=\alpha+\lambda<\sup(\card{M}).
\]
So, $\delta\le \sup(\card{M})$, and the other side of the inequality
is shown similarly.]

{\noindent (2)} If $\bar{M}=\langle M_i:\,i<i^\ast\rangle$ is a
$\le_{K_{\rm ap}}$-increasing sequence, and $\lambda |\delta$, then
$\langle M_i\rest\delta:\,i<i^\ast\rangle$ is $\le_{K_{\rm ap}}$-increasing,
by Definition \ref{general}(1)(c), and if $i^\ast<\lambda$,
\[
 \bigcup_{i<i^\ast}(M_i\rest\delta)=
(\bigcup_{i<i^\ast}M_i)\rest\delta
\]
is the $\le_{K_{\ap}}$-lub of
$\langle M_i\rest\delta:\,i<i^\ast\rangle$, by (3)(d) in Definition \ref{general}.

{\noindent (3)} Suppose $M_l$ for $l<3$ are as in Definition \ref{general}(3)(i)
(amalgamation).
Then we can without loss of generality assume that $M\rest \Ev=M_1\rest \Ev$, as clearly
there is a lawful mapping $g:\,M\into M^\ast$
extending ${\rm id}_{M_1}$ for some $M^\ast$ with
$M^\ast\rest \Ev=M_1\rest \Ev$. 

{\noindent (4)} Suppose that $M_0, M_1$ and $M_2$ are as in Definition
\ref{general}(3)(e) (end extension/amalgamation). Then we can assume
$M_3\subseteq \iota(M_2)$, as by Definition \ref{general}(1)(c) we can replace
$M_3$ by $M_3\rest \iota(M_2)$.
\end{Remark}

\begin{Notation} Suppose that $K_{\ap}$
is a weak $\lambda$-approximation family and 
$\bar{\tau}$ is a sequence of vocabularies as in Definition
\ref{general}(1)(a). We say that $K_{\ap}$ {\em is written in}
$\bar{\tau}$.
\end{Notation}

\begin{Definition}\label{md} [Sh 457] 
\begin{description}
\item{(1)} Let $(K_{\ap},\le_{K_{\ap}})$
be a
weak $\lambda$-approximation family and $\Gamma\subseteq K_{\ap}$.
We say that $\Gamma$ is $(<\lambda)$-{\em closed} iff
for every $\le_{K_{\ap}}$-increasing 
chain of size $<\lambda$ of elements of $\Gamma$, the lub of the chain is in $\Gamma$.

\item{(2)} Suppose that $(K_{\ap},\le_{K_{\ap}})$ is a
weak $\lambda$-approximation family. We let
\[
K_{{\rm md}}^-=K_{{\rm md}}^-[K_{\ap}]
\deq\left\{\Gamma:\,\begin{array}{l}
{\rm (i)}\,\Gamma\mbox{ is a }(<\lambda)
\mbox{-closed subset of }K_{\ap},\\
{\rm (ii)}\,\Gamma\mbox{ is }\le_{K_{\ap}}\mbox{-directed},\\
{\rm (iii)}\,\mbox{ for cofinally many } \beta<\lambda^+\mbox{ we have}\\
(\exists M\in \Gamma)(\exists
\gamma\in\card{M}) \iota(\gamma)=\iota(\beta)\\
\quad (\mbox{e.g.} \gamma=\beta)
\end{array}
\right\}.
\]
We let
$K_{{\rm md}}=K_{{\rm md}}[K_{\ap}]
\deq$
\[
\left\{\Gamma\in K_{{\rm md}}^-:\,
\begin{array}{l}
{\rm (iv)}\, (M\in\Gamma\,\,\&\,\,M\le_{K_{\ap}} M_1)\implies\\
\,\,\,(\exists M_2\in \Gamma)
(\exists h\mbox{ lawful})[h:\,M_1\to M_2\\
\quad\quad\quad\mbox{ embedding over }M]\\
{\rm (v)}\,M\in \Gamma\,\,\&\,\,N\le M\implies N\in \Gamma
\end{array}
\right\}.
\]

\item{(3)} If $K_{\ap}$ is as above and $\alpha<\lambda^+$, we define
$K^-_{\md}[K_{\ap}^\alpha]$ as the set of $\Gamma\subseteq K_{\ap}$
such that
\begin{description}
\item{(a)} $M\in\Gamma\implies\card{M}\subseteq\alpha$,
\item{(b)} $\Gamma$ satisfies (i) --(ii) from (2) above.
\end{description}
Similarly for $K_{\md}[K_{\ap}^\alpha]$.
\end{description}
\end{Definition}

\begin{Claim}\label{10B} Suppose that $\Gamma\in K_{\rm md}[K_{\rm ap}]$,
while $N\in \Gamma$ and $h:\,N\into M$ is a lawful embedding.
Then there is $N'\in \Gamma$ and a lawful embedding $g:\,M
\into N'$ such that for $x\in N$ we have $g(h(x))=x$.
\end{Claim}

\begin{Proof of the Claim} There is a lawful isomorphism $f:\,M\into M'$
for some $M'\ge N$ such that $f(h(x))=x$ for all $x\in N$. Then by
(iv) in the definition of $K_{\rm md}$, there is $N'\in \Gamma$ and a lawful
embedding $g':\,M'\into N'$ such that $g'\rest N={\rm id}_N$.

Let $g:\,M\into N'$ be given by letting $g(x)=g'(f(x))$, so $g$ is
a lawful embedding and for $x\in N$ we have $g(h(x))=g'(f(h(x))=g'(x)=x$.
$\eop_{\ref{10B}}$
\end{Proof of the Claim}
\begin{Definition}\label{absel}
(1) $\KK=(K,\le_{\KK})$ is an {\em abstract elementary class} iff
$K$ is a class of
models of some fixed vocabulary $\tau=\tau_{\KK}$ and $\le_{\KK}=
\le_K$
is a two place relation on $K$, satisfying the following axioms:
\begin{description}
\item{Ax 0:} If $M\in K$, then all $\tau$-models isomorphic to $M$
are also in $K$. The relation $\le_K$ is preserved under isomorphisms,
\item{Ax I:} If $M\le_K N$, then $M$ is a submodel of $N$,
\item{Ax II:} $\le_K$ is a partial order on $K$,
\item{Ax III, IV:} The union of
a $\le_K$-increasing continuous chain $\bar{M}$ of elements of $K$ is
an element of $K$, and the lub of $\bar{M}$ under $\le_K$,
\item{Ax V:} If $M_l\le_K N$ for
$l\in\{0,1\}$ and $M_0$ is a submodel of $M_1$, then $M_0\le_K M_1$,
\item{Ax VI:}
There is a cardinal $\kappa$ such that for every $M\in\KK$ and $A
\subseteq \card{M}$, there is $N\le_K M$
such that
$A\subseteq\card{N}$ and
$\card{\card{N}}\le \kappa\cdot(\card{A}+1)$. The least
such $\kappa$ is denoted by
${\rm LS}(\KK)$ and called the L\"owenheim-Skolem number of $\KK$.
\end{description}

{\noindent (2) If $\lambda$ is a cardinal and $\KK$ an abstract
elementary class,
we denote by $\KK_\lambda$ the family of all elements of $\KK$
whose cardinality is $\lambda$.}

{\noindent (3)} For $\KK$ an abstract elementary class, and
$\lambda$ a cardinal, we say that $\KK_\lambda$ {\em
has a universal} iff there is $M^\ast\in \KK_\lambda$
such that for all $M\in \KK_\lambda$ we have
that some $M'$ which is isomorphic to $M$ satisfies $M'\le_{\KK}
M^\ast$. Such $M^\ast$ is called {\em universal for} $\KK_\lambda$.

{\noindent (4)} Suppose that $\KK$ is an abstract elementary class. We shall say that
a member $M$ of $\KK$ is $\le_{\KK}$-{\em embeddable} in a member $N$ of $\KK$
iff there is an isomorphism between $M$ and some $M'\in\KK$ satisfying
$M'\le_{\KK} N$.
\begin{description}
\item{(a)} $\KK$ is said to have the {\em joint embedding property}
iff for any $M_1, M_2\in \KK$, there is $N\in \KK$ such that
$M_1,M_2$ are
$\le_{\KK}$-embeddable into $N$.
\item{(b)} $\KK$ is said to have {\em amalgamation} iff for all $M_0, M_1, M_2
\in \KK$ and $\le_{\KK}$-embeddings $g_l:\,M_0\into M_l$ for $l\in\{1,2\}$,
there is $N\in \KK$ and $\le_{\KK}$-embeddings $f_l:\,M_l\into N$ such that
$f_1\circ g_1=f_2\circ g_2$.
\end{description}
Similar definitions are made to describe when $\KK_\lambda$ has the joint embedding
property or amalgamation.
\end{Definition}

\begin{Convention} We shall only work with abstract elementary classes
which have the joint embedding property and amalgamation.
\end{Convention}

\begin{Note}\label{duga}
The following notes are not hard and the proofs are to be found
in \cite{Sh88}. We include them here for the reader's convenience.

{\noindent (1)} Suppose that $\KK$ is an abstract elementary class.
If
$\bar{M}=\langle M_i:\,i<\delta\rangle$ is a
$\le_{\KK}$-increasing chain (not necessarily continuous), then 
$\bigcup_{i<\delta}M_i$ is the $\le_{\KK}$-lub of $\bar{M}$.

[Why? Prove this by induction on $\delta$. The nontrivial case is when
$\delta$ is a limit. Define for $i<\delta$ a model $N_i$ to be $M_i$
if $i$ is non-limit, and $\bigcup_{j<i}M_j$ otherwise. Now
$\bar{N}=\langle N_i:\,i<\delta\rangle$ is increasing continuous and
$\bigcup_{i<\delta}N_i=\bigcup_{i<\delta}M_i$ is the lub of $\bar{N}$,
hence of $\bar{M}$.]

{\noindent (2)} If $\KK$ is an abstract elementary class, $\KK$ is closed
under unions of $\le_{\KK}$-directed subsets, and the union
of a $\le_{\KK}$-directed
subset of $\KK$ is the $\le_{\KK}$-lub of it.

[Why? By induction on $\kappa$, we prove that for any $\DD\subseteq\KK$
which is $\le_{\KK}$-directed and has size $\kappa$, the $\le_{\KK}$-lub
of $\DD$ is $\bigcup\DD$. For $\kappa\le\aleph_0$, this is clear. If $\kappa$
is a limit $>\aleph_0$, 
let $\langle\kappa_\alpha:\,\alpha<\cf(\kappa)\rangle$ be cofinal
increasing to $\kappa$, each $\kappa_\alpha$ regular, and $\DD=\bigcup_{\alpha<
\cf(\kappa)}
\DD_\alpha$, where each $\DD_\alpha$ is $\le_{\KK}$-directed
and has size $\kappa_\alpha$, and $\DD_\alpha$'s
are $\subseteq$-increasing. Now apply the induction hypothesis and (1). If
$\kappa=\lambda^+$, then we can find $\langle \DD_\alpha:\,\alpha<\lambda^+
\rangle$ increasing to $\DD$, each $\le_{\KK}$-directed and of size
$\le\lambda$.]
\end{Note}

\begin{Definition}\label{tendsto} Suppose that $\KK$ is an abstract elementary class
with $\tau_\KK=\tau$,
and $K_{\ap}$ is a weak [strong] $\lambda$-approximation family
written in
\[
\langle \tau_i:\,i<\lambda^+\,\,\&\,\,\lambda|i\rangle,
\]
such that 
\begin{description}
\item{(1)} For all $i$, we have $\tau\subseteq \tau_i$,
\item{(2)} $M\in K_{\ap}\implies M\rest\tau\in \KK$,
\item{(3)} $M\le_{K_{\ap}}N\implies M\rest\tau\le_{\KK} N\rest\tau$ and
\item{(4)} For every $M\in\KK$ with $\norm{M}<\lambda$ 
there is $N\in K_{\ap}$ such that
\[
M\mbox{ is }\le_{\KK}\mbox{-embeddable into } N\rest\tau.
\]
\end{description}
We say that $K_{\ap}$ {\em tends to [strongly] $\lambda$-approximate} $\KK$.

We may just say ``$K_{\ap}$ tends to approximate $\KK$" if the rest is clear
from the context.
\end{Definition}

\begin{Observation}\label{union} Suppose that
$K_{\ap}$ is a strong
$\lambda$-approximation family
which tends to approximate $\KK$
and $\Gamma\in K^-_{{\rm md}}$.
Then
\begin{description}
\item{(1)}
$M_\Gamma$ defined by letting
\[
M_\Gamma\deq \bigcup_{M\in \Gamma}M\rest\tau
\]
is an element of $\KK$
and for every $M\in \Gamma$ we have
$M\rest\tau\le_{\KK} M_\Gamma$,
and in fact $M_\Gamma$ is the $\le_{\KK}$-lub of
$\{M\rest\tau:\,M\in\Gamma\}$.

\item{(2)} For every $\Gamma,\Gamma^\ast\in K_{\rm md}^-[K_{\rm ap}]$ such
that $\Gamma\subseteq\Gamma^\ast$, we have $M_\Gamma\le_{\KK}M_{\Gamma^\ast}$.
\end{description}

[Why? (1)
As $\{M\rest\tau:\,M\in\Gamma\}$ is $\le_{\KK}$-directed.

(2) By (1) and Note \ref{duga}(2).]
\end{Observation}

\begin{Notation}\label{grande} Suppose that
an approximation family $K_{\ap}$
tends to approximate $\KK$,
while $\Gamma\in K^-_{{\rm md}}$. If we
write $M_\Gamma$, we always mean the
model obtained from $\Gamma$ as in Observation \ref{union}.

\end{Notation}

\begin{Definition}\label{priblizno}
Let $K_{\ap}$ be a strong $\lambda$-approximation
family which tends to $\lambda$-approximate $\KK$ and let $\KK^+$ be a subclass of
$\KK_{\lambda^+}$. Assume
\begin{description}
\item{($\ast$)} For every $M^\ast\in\KK^+$, there is $\Gamma\in
K^-_{\md}[K_{\ap}]$ with $\{\card{M}:\,M\in\Gamma\}$ a club of
$[\Ev]^{<\lambda}$ such that for some $M'$ isomorphic to
$M^\ast$, we have $M'\le_{\KK}M_\Gamma$.
\end{description}
Then we say that $K_{\ap}$ {\em approximates} $\KK^+$.

\end{Definition}

\begin{Claim}\label{tree} Suppose that
\begin{description}
\item{(1)} $\lambda\le\kappa$,
\item{(2)}
$\KK$ is an abstract elementary
class,
\item{(3)} LS$(\KK)\le\kappa$ and $\KK_\kappa$ has amalgamation,
\item{(4)}
$\TT\subseteq{}^{<\lambda^+}(\lambda^+)$ ordered by $\unlhd$
(i.e. being an initial segment) is a tree with each level of size $\le\lambda^+$,
\item{(5)} For $\eta\in \TT$ we have $M_\eta\in \KK$, so that
\[
\eta\unlhd\nu\implies M_\eta\le_{\KK}M_\nu,
\]
\item{(6)} $\eta\in\TT\implies \norm{M_\eta}=\kappa$.
\end{description}
Then there are $M^\ast=M^\ast[\TT]\in \KK$
and $\langle g_\eta:\,\eta\in \TT\rangle$
such that
\begin{description}
\item{(A)} For all $\eta\in \TT$ we have that $g_\eta$ is a $\KK$-embedding
from
$M_\eta$ into $M^\ast$,
\item{(B)} $\eta\le\nu\implies g_\eta\subseteq g_\nu$,
\item{(C)} $\norm{M^\ast}\le\kappa\cdot\lambda^+$.
\end{description}
({\scriptsize{The intended use of this claim is when $\kappa=\lambda$.}})
\end{Claim}

\begin{Proof of the Claim} For $i^\ast\le\lambda^+$, let $\TT\rest i^\ast\deq
\TT\cap {}^{<i^\ast}\lambda^+$.

By induction on $i^\ast\le$ the height of $\TT$, we prove that
$M^\ast[\TT\rest i^\ast]$ and
$\langle g_\eta^{i^\ast}:\,\eta\in \TT\rest i^\ast\rangle$
can be defined to satisfy (A)-(C) with $\TT\rest i^\ast$ in place of $\TT$
and $M^\ast[\TT\rest i]$ in place of $M^\ast$, and so that 
\[
i\le i^\ast\implies M^\ast[\TT\rest i]\le_{\KK}M^\ast[\TT\rest i^\ast],
\]
\[
[\eta\in \TT\rest i\,\,\&\,\,i\le i^\ast]\implies g_\eta^i=g_\eta^{i^\ast}.
\]

\underline{$i^\ast=0$}. Trivial.

\underline{$i^\ast=i+1$}. Let ${\rm lev}_i(\TT)=
\{\eta_j:\,j< j^\ast\le\lambda^+\}$. For simplicity in notation we
assume that $j^\ast$ is a limit > 0, the other cases are similar.
By induction on $j$ we build
$\langle M^\ast_j:
\,j<j^\ast\rangle$ so that $M^\ast_0=M^\ast[\TT\rest i]$ and $M_{\eta_j},
M^\ast_j\le_{\KK}M^\ast_{j+1}$, while $\norm{M^\ast_j}=\kappa$,
and $M^\ast_\delta=\bigcup_{j<\delta}M^\ast_j$ for
$\delta$ a limit. We use 
amalgamation and the induction hypothesis to obtain (B). Namely,
to define $M^\ast_{j+1}$, let first $M'_j\deq\bigcup\{M_\nu:\,\nu
\initial\eta_j\}$ and $g_j\deq\cup\{g_\nu^i:\,\nu
\initial\eta_j\}$, which is well defined by the induction hypothesis.
Hence $g_j:\,M'_j\into M_0^\ast\le_{\KK}M^\ast_j$ is a $\le_{\KK}$-embedding,
as is ${\rm id}:\,M'_j\into M_{\eta_j}$. Using amalgamation,
we can find $M^\ast_{j+1}\in \KK$ and $\le_{\KK}$-embeddings
$f:\,M_j^\ast\into M^\ast_{j+1}$ and $g_{\eta_j}:\,M_{\eta_j}\into M^\ast_{j+1}$
such that $f\circ g_j=g_{\eta_j}\rest M'_j$. By Ax 0 of Definition \ref{absel},
without loss of generality we have $f={\rm id}$. By Ax VI of the same
Definition, we can also assume that $\norm{M^\ast_{j+1}}\le\kappa$.
Now let $M^\ast[\TT\rest i^\ast]\deq\bigcup_{j<j^\ast}M^\ast_j$.

\underline{$i^\ast$ a limit}. $M^\ast[\TT]=\bigcup_{i<i^\ast}M^\ast[\TT\rest i]$.

$\eop_{\ref{tree}}$
\end{Proof of the Claim}

\begin{Observation}\label{dodatna} With the notation of Claim \ref{tree}, if $\rho$
is a branch of $\TT$, then $\bar{M}=\langle M_\eta:\,\eta\in\TT
\,\,\&\,\,\eta\initialeq\rho\rangle$ is a $\le_{\KK}$-increasing chain of
$\KK$. Hence $\bigcup \bar{M}$ is the $\le_{\KK}$-lub of $\bar{M}$, and so
$\bigcup \bar{M}$ is $\le_{\KK}$-embeddable into $M^\ast[\TT]$.
\end{Observation}

\begin{Definition}\label{simplicity} For a strong $\lambda$-approximation
family $K_{\ap}$ we say that it is
{\em workable}
iff
for every
$\Gamma\in K_{{\rm md}}^-[K_{\ap}]$ such that
$M\in \Gamma\implies\card{M}\subseteq \Ev$,
for all
$\delta_1<\delta_2\in S^{\lambda^+}_\lambda$ the following holds:

Suppose that for $l\in \{1,2\}$ we are given $(M_l, N_l)$
such that
\begin{description}
\item{(i)} $M_l\in \Gamma$,
\item{(ii)} $M_l\le_{K_{\ap}} N_l\in K_{\ap}$, 
\item{(iii)} $\card{N_l}\cap\{2\beta:\,\beta<\lambda^+\}=\card{M_l}$,
\item{(iv)} $\card{N_1}\subseteq \delta_2$,
\item{(v)} $N_1\rest\delta_1=N_2\rest\delta_2$,
\item{(vi)} Some $h$ is a
lawful $K_{\ap}$-isomorphism from $N_1\rest\tau(\delta_1)$
onto $N_2\rest\tau(\delta_1)$ mapping $M_1$ onto $M_2$,
\item{(vii)} $h\rest(N_1\rest\delta_1)$ is the identity.
\end{description}
Then
there are $M\in \Gamma$
and $N\in K_{\ap}$ with $M\le N$, and $g_l$
for $l\in \{1,2\}$ such that $M_l\le M\le N$ and
$g_l$ is a $\le_{K_{\ap}}$-embedding of $N_l$ into $N$,
with $g_l\rest M_l={\rm id}_{M_l}$. In addition,
$\card{N}\cap\{2\beta:\,\beta<\lambda^+\}=\card{M}$ and $g_l\rest(N_l\rest\delta_l)$
is fixed.
\end{Definition}

\begin{Note}
For those familiar
with definitions in [Sh 457],
we emphasize that smoothness was assumed throughout. That is, our
definition of $K_{\rm ap}$ is less general than the one in \cite{Sh457}, and any
strong $\lambda$-approximation family in the sense of our Definition 
\ref{md} automatically satisfies the condition which in \cite{Sh457} was called
smoothness.
\end{Note}

\section{Universals in $\lambda^+$.}\label{forcing}

\begin{Definition}\label{thegame} [Sh 546] Suppose that
$\lambda>\aleph_0$ is a cardinal and $\varepsilon<\lambda$ a limit
ordinal.
A forcing notion $Q$ satisfies $\ast^\varepsilon_\lambda$ iff
player I has a winning strategy in the following game
$\ast^\varepsilon_\lambda[Q]$:

\underline{Moves:} The play lasts $\varepsilon$ moves. For
$\zeta<\varepsilon$, the $\zeta$-th move is described by:
\begin{description}
\item{\underline{Player I}}: If $\zeta\neq 0$, I chooses $\langle
q_i^\zeta:\,
i<\lambda^+\rangle$ such that $q_i^\zeta\in Q$
and $q^\zeta_i\ge p^\xi_i$ for all $\xi<\zeta$, as well as a
function $f_\zeta:\,\lambda^+\to\lambda^+$ which is regressive on
$C_\zeta\cap S^{\lambda^+}_\lambda$ for some
club $C_\zeta$ of $\lambda^+$.
If $\zeta=0$, we let $q_i^\zeta\deq\emptyset_Q$ and $f_\zeta$
be identically 0.
\item{\underline{Player II}}: Chooses $\langle p_i^\zeta:\,i<\lambda^+
\rangle$
such that $q_i^\zeta\le p_i^\zeta\in Q$ for all $i<\lambda^+$.
\end{description}

\underline{The Outcome:} Player I wins iff:

For some club $E$
of $\lambda^+$, for any $i<j\in E\cap S^{\lambda^+}_\lambda$,
\[
\bigwedge_{\zeta<\varepsilon}f_\zeta(i)=f_\zeta(j)
\implies [\{p^\zeta_i:\,\zeta<\varepsilon\}\cup
\{p^\zeta_j:\,\zeta<\varepsilon\}\mbox{ has an upper bound in }Q].
\]
We say that $E\subseteq\bigcap_{\zeta<\varepsilon}C_\zeta$ is a witness that I
won.


{\noindent (2)} A {\em winning strategy for I in }$\ast^\varepsilon_\lambda[Q]$
is a function St$
=({\rm St}_\ast, {\rm St}^\ast)$ such that in any play
\[
\left\langle\langle q^\zeta_i:\,i<\lambda^+\rangle, f_\zeta,
\langle p^\zeta_i:\,i<\lambda^+\rangle:\,\zeta<\varepsilon\right\rangle
\]
in which we have for all $\zeta,i$
\[
q^\zeta_i={\rm St}_\ast(i,\left\langle 
\langle p^\xi_j:\,j<\lambda^+\rangle:\,\xi<\zeta\right\rangle),
f_\zeta={\rm St}^\ast(\left\langle
\langle p^\xi_j:\,j<\lambda^+\rangle:\,\xi<\zeta\right\rangle),
\]
I wins.

{\scriptsize{i.e. a winning strategy for I depends only on
the moves of II and $f_\zeta$ and $C_\zeta$ can be defined from 
$\langle\langle p^\xi_j:\,j<\lambda^+\rangle:\,\xi<\zeta\rangle$.}}
\end{Definition}

\begin{Fact}\label{unnecessary} [Sh 546] Suppose that
$\lambda>\aleph_0$ is
a cardinal
satisfying $\lambda^{<\lambda}=\lambda$, and $\varepsilon<\lambda$ a limit ordinal.

{\noindent (1)} If $P$ is a forcing notion satisfying
$\ast^\varepsilon_\lambda$, then $P$ satisfies $\lambda^+$-cc.

{\noindent (2)} Suppose that $P$ is the result of an iteration of
$(<\lambda)$-complete forcing satisfying $\ast^\varepsilon_\lambda$.
Then $P$ is $(<\lambda)$-complete and satisfies 
$\ast^\varepsilon_\lambda$.
\end{Fact}

\begin{Proof of the Fact} (1) Suppose that $\bar{p}=\langle p_i:\,i<\lambda^+
\rangle$ is a sequence of elements of $P$,
and consider a game of $\ast^\varepsilon_\lambda[P]$ in which
II plays $\bar{p}$ as the first move, and I plays according to a
winning strategy. At the end of the game, let $E$ be
a club of $\lambda^+$ witnessing that I won, and let $i<j$
be in $E\cap S^{\lambda^+}_\lambda$ such that for all
$\zeta<\varepsilon$ we have that $f_\zeta(i)=f_\zeta(j)$,
which exists as these functions are regressive. We in particular
obtain that $p_i$ and $p_j$ are compatible in $P$.

{\noindent (2)} We refer the reader to [Sh 546].
\end{Proof of the Fact}

\begin{Theorem}\label{universal} Suppose that the following are
satisfied in a universe $V_0$ of set theory:
\begin{description}
\item{(A)} $\aleph_0\le \lambda=\lambda^{<\lambda}<\lambda^+=2^\lambda<
2^{\lambda^+}\le \kappa<\mu=
\mu^{\kappa}$,
\item{(B)} $R^\ast$ is the forcing notion which adds $\mu$ many Cohen subsets 
$\langle \rho^\ast_\alpha:\,\alpha<\mu\rangle$ to $\lambda^+$ by conditions
of size $\le \lambda$.

\item{(C)} $\TT={}^{<\lambda^+}(\lambda^+)$ of $V_0$, ordered by ``being
an
initial segment",
\item{(D)} If $\lambda>\aleph_0$, we are given a limit
ordinal $\varepsilon<\lambda$. \end{description}

Then in $V\deq V^{R^\ast}_0$ for some $P$ we have
\begin{description}
\item{(a)} $P$ is a forcing notion of cardinality $\mu$,
\item{(b)} $P$ is $(<\lambda)$-complete and $\lambda^+$-cc
(and if $\lambda>\aleph_0$, $P$ satisfies $\ast^\varepsilon_\lambda$),
\item{(c)} In $V^P$ we have $\lambda^{<\lambda}=\lambda$ and $2^\lambda=
2^{\lambda^+}=\mu$,
\item{(d${}_0$)} If $\lambda=\aleph_0$, then $MA(\aleph_1)$ holds in $V^P$,
\item{(d${}_1$)} If $\lambda>\aleph_0$, then the following holds in
$V^P$: if $Q$ is a $(<\lambda)$-complete forcing notion of cardinality
$<\kappa$ and
satisfies $\ast^\varepsilon_\lambda$, and if
we are given a family $\{\II_j:\,j<\lambda^+\}$
of dense subsets of $Q$, then for some directed $G
\subseteq Q$ we have that $G\cap \II_j\neq\emptyset$ for all
$j<\lambda^+$,
\item{(e)} In $V^P$, if $K=K_{\ap}$ is a workable
strong $\lambda$-approximation
family, then
we can find 
\[
\langle \bar{\Delta}_\beta=\langle
\Delta_\eta^\beta:\,\eta\in\TT\rangle:\,\beta<\lambda^{++}\rangle
\]
such that
\begin{description}
\item{(i)} For every $\beta<\lambda^{++}$ and  $\eta\in\TT$ we have
$\Delta_\eta^\beta\subseteq K_{\ap}^{\lambda\cdot\lg(\eta)}$ is
$\le_{\KK_{\ap}}$-directed, and also for $\eta\initialeq\nu\in \TT$, we have that
\[
\Delta_\eta^\beta=\{M\rest(\lambda\cdot\lg(\eta)):\,M\in \Delta_\nu^\beta\},
\]
\item{(ii)} For any $\lambda^+$-branch $\rho$ of $\TT$ and
$\beta<\lambda^{++}$, we have
\[
\bigcup\{\Delta_\eta^\beta:\,\eta\unlhd\rho\}\in K^-_{\md}[K_{\ap}],
\]
\item{(iii)} For any $\Gamma\in K_{\md}[K_{\rm ap}]$, for some
$\beta<\lambda^{++}$
we have that
$M_\Gamma$ is isomorphically embeddable into
$M_{\bigcup_{i<\lambda^+}\Delta^\beta_{\rho\rest i}}$
for some $\lambda^+$-branch $\rho$ of $\TT$ with $\rho\in V$ (for the notation
see \ref{grande}), \end{description}
\item{(f)}
In $V^P$, if $K_{\ap}$ is a workable
strong $\lambda$-approximation
family and $\Gamma^-$
is an element of $K^-_{\md}[K_{\ap}]$
such that
$M\in \Gamma^-\implies\card{M}\subseteq
\Ev$,
then there is $\Gamma\in K_{\md}[K_{\ap}]$ such that
$\Gamma^-\subseteq\Gamma$.
\end{description}
\end{Theorem}

Once we prove the theorem, we shall be able to draw the following

\begin{Conclusion}\label{conclus} Suppose that $V$ satisfies
\[
\aleph_0\le\lambda=\lambda^{<\lambda}<\lambda^+=2^\lambda<
2^{\lambda^+}\le \kappa<\mu=
\mu^{\kappa},
\]
and if $\lambda>\aleph_0$, we are given a limit ordinal $\varepsilon<
\lambda$.
Then there is a 
cofinality and cardinality preserving forcing extension $V^\ast$ of $V$ which
satisfies
\begin{description}
\item{(1)} For every abstract elementary class $\KK$ for
which there is a workable
$\lambda$-approximation family $K_{\ap}$ which approximates $\KK$,
and such that ${\rm LS}(\KK)\le\lambda$,
there are $\lambda^{++}$ elements
$\{M_\alpha:\,\alpha<\lambda^{++}\}$ of $\KK_{\lambda^+}$ which are
jointly
universal for ${\KK}_{\lambda^+}$,

\item{(2)} $\aleph_0\le \lambda^{<\lambda}=\lambda<2^\lambda=2^{\lambda^+}=\mu=
\mu^{\kappa}$,

\item{(3)(a)} In the case $\lambda=\aleph_0$: $MA(\aleph_1)$ holds,

\item{(3)(b)} In the case $\lambda>\aleph_0$: if $Q$ is a
$(<\lambda)$-complete forcing notion of cardinality
$<\kappa$,
satisfying $\ast^\varepsilon_\lambda$, and
we are given a family $\{\II_j:\,j<\lambda^+\}$ of
dense subsets of $Q$, then for some directed $G
\subseteq Q$ we have that $G\cap \II_j\neq\emptyset$ for all
$j<\lambda^+$,

\item{(4)} If $\KK$ is an abstract elementary class with
${\rm LS}(\KK)\le\lambda$ and $\KK^+$ is a subclass
of $\KK_{\lambda^+}$ for which there is a workable strong
$\lambda$-approximation family $K_{\ap}$ which approximates $\KK^+$,
and such that for every tree $\TT$ of the form from
Claim \ref{tree} in which every $M_\eta$ is the union of
$\le\lambda$ elements of $K_{\ap}$ we have that $M^\ast[\TT]\in\KK^+$,
then 
there are $\lambda^{++}$ elements
$\{M_\alpha:\,\alpha<\lambda^{++}\}$ of $\KK^+$ which are
jointly
universal for $\KK^+$.
\end{description}
\end{Conclusion}

\begin{Remark} The informal plan of the proof of the theorem
and the conclusion is as follows. The
purpose of forcing with $R^\ast$ is to
make $2^{\lambda^+}=\mu$ and add $\mu$ branches through
$\TT$. Then $P$ will be an iteration
of $\lambda^{++}$ blocks of $\mu$ steps each.
Hence
\[
P=\langle
P_\alpha,\name{Q}_\beta:\,\alpha\le\lambda^{++},\beta<\lambda^{++}\rangle
\]
and for each $\beta$ we have $\name{Q}_\beta=\langle
Q^\beta_i,\name{R}^\beta_j:\,i
\le\mu, j<\mu\rangle$. Each $\name{R}^\beta_j$ will be one of four possible
kinds (three in case $\lambda=\aleph_0$). Let us first describe the situation 
when $\lambda>\aleph_0$.

At kind 1 coordinates we shall be taking care of
the form of Martin's Axiom given in (d) of the Theorem. Each kind 2 coordinate
$\name{R}^{\beta+1}_j$ takes a workable strong $\lambda$-approximation family
$K_{\rm ap}$ from $V^{P_\beta}$ and a family of $\le \mu$ elements of $K_{\rm
ap}[K_{\rm md}]$
and introduces a tree of elements of $K_{\rm ap}$ indexed by $\TT$, which gives
$\bar{\Delta}_\beta$ as in (e)(i)-(ii) of the Theorem. This tree will
also have the
property that for every $\Gamma\in K_{\rm
ap}[K_{\rm md}]\cap V^{P_\beta}$, there is a branch $\rho$ of $\TT$ with
$\rho\in V$ and a tree $T$ whose elements are pairs $(N,h)$ with $N\in \Gamma$
and $h$ an embedding from $N$ into
$M_{\bigcup_{i<\lambda^+}}\Delta^\beta_{\rho\rest i}$, ordered by extension.
Then for some $\beta'\in(\beta,
\lambda^{++})$, there will be a forcing of the third kind that will introduce a
branch through $T$ and so have (e)(iii) of the Theorem.

At the remaining coordinates, for a
workable strong $\lambda$-approximation family $K$ introduced at some
earlier stage,
we embed $M_{\Gamma^-}$ for some $\Gamma^-\in K_{\md}^-$ into
$M_{\Gamma}$ for some $\Gamma\in K_{\md}$.

If $\lambda=\aleph_0$, the forcing is easier because we do not need a strong
chain condition in order to be able to iterate. So the kind two coordinates,
which satisfy ccc but not the stronger analogue of it needed
if $\lambda>\aleph_0$, are simplified and guarantee (e)(i)-(iii) immediately.
This eliminates the need for kind three coordinates.

To get the conclusion for a given $\KK$ as in (1), recall
from \S\ref{definitions} that for every $\beta<\lambda^{++}$  there is a model
$M^\ast_\beta$ in $\KK_{\lambda^+}$ such that for every branch $\rho$ of $\TT$,
the model read along the branch in the tree indexed by $\bar{\Delta}_\beta$,
embeds into $M^\ast_\beta$. As $K_{\ap}$ approximates $\KK$ (see Definition
\ref{priblizno}), for every $M\in\KK$, there is $\Gamma^-\in K_{\md}^-$ (and of
the kind required by the Theorem) such that $M$ embeds into $M_{\Gamma^-}$. From the Theorem,
there is $\Gamma\in K_{\md}$ such that $M_{\Gamma^-}$ embeds into $M_{\Gamma}$.
This $\Gamma$ is in $V^{P_\beta}$ for some $\beta<\lambda^{++}$, and hence some 
$\name{R}^{\beta+1}_j$ will guarantee that
$M_\Gamma$ embeds into $M^\ast_\beta$.

The {\em proof of the Conclusion} is given after the proof
of the Theorem,
close to the end of the section.
\end{Remark}

\begin{Proof of the Theorem} Let $R^\ast$ be as in the statement of the
Theorem, and let $V=V^{R^\ast}$. Then in $V$ we clearly have
$\aleph_0\le \lambda=\lambda^{<\lambda}$, while $2^\lambda=\lambda^+$,
$2^{\lambda^+}=\mu$ and the cardinalities and cofinalities of $V_0$ are
preserved. 

Let $\langle f^\ast_\alpha:\,\alpha<\mu\rangle$ list the $\lambda^+$-branches
of $\TT$ in $V$.

We make some easy observations:

\begin{Note}\label{dominating} (1)
It suffices to prove the conclusion weakened by
requiring each $Q$ being considered in (d)${}_1$, to have the set of elements
some ordinal $<\kappa$.  

{\noindent (2)} By renaming, each $K_{\rm ap}$ considered in the theorem
can be assumed to have its vocabulary included in $\HH(\lambda^+)$.

\end{Note}

\begin{Definition}\label{candidate}
We define $P$ as $P_{\lambda^{++}}$ in the iteration
\[
\bar{P}=\langle P_\alpha,\name{Q}_\beta:\,
\alpha\le\lambda^{++},\beta<\lambda^{++}\rangle,
\]
where 
\begin{description}
\item{$(\alpha)$} $\bar{P}$ is a $(<\lambda)$-support iteration.
\item{$(\beta)$} For each $\beta<\lambda^{++}$, in $V^{P_\beta}$ we have that
$Q_\beta$ is $Q^\beta_\mu$ in the iteration
\[
\bar{Q}^\beta=\langle Q^\beta_i,\name{R}^\beta_j:\,
i\le\mu ,j<\mu\rangle,
\]
where:
\begin{description}
\item{(i)} the iteration in $\bar{Q}^\beta$ is made with $(<\lambda)$-supports,
\item{(ii)} 
for each $j<\mu$ one of the following
occurs:
\end{description}
\smallskip

\underline{Case 1}. $\name{R}^\beta_j$ is a $Q^\beta_j$-name of
a $(<\lambda)$-complete forcing notion
which satisfies $\ast^\varepsilon_\lambda$ if $\lambda>\aleph_0$,
and is ccc if $\lambda=\aleph_0$; and whose
set of elements is some ordinal
$<\kappa$.

\underline{Case 2}. For some $P_\beta$-name of
a workable
$\lambda$-approximation family $\name{K}^\beta_{{\rm ap},j}$,
abbreviated as $\name{K}$, and
elements
$\{\Gamma_\alpha=\Gamma^{\beta,j}_\alpha:\,\alpha<\mu\}$ of $K_{\rm md}^{V^{P_\beta}}$,
we have that $\name{R}=\name{R}^\beta_j$ is defined as follows. We work in
$V^{P_\beta\ast\name{Q}^\beta_j}$. For $M\in K$ we let
\[
w[M]\deq\{\gamma,\gamma+1:\,M\cap[\lambda\gamma,\lambda(\gamma+1))
\neq\emptyset\}.
\]

\underline{Subcase 2A}. $\lambda=\aleph_0$. The elements of $R$ are conditions
of the form
\[
p=\langle u^p,\langle M^p_\eta:\,\eta\in u^p\rangle, b^p,
\langle c^p_\alpha:\,\alpha\in b^p\rangle,
\langle (N^{p}_{\alpha,\iota}, h^p_{\alpha,\iota}):\,\alpha\in
b^p,\iota\in c^p_\alpha\rangle\rangle,
\]
where
\begin{description}
\item{(a)}[{\em closure under intersections}] $u=u^p\in [\TT]^{<\lambda}$ is
closed under intersections
\item{(b)} $\eta\in u\implies M^p_\eta\in K\,\,\&\,\,\card{M^p_\eta}
\subseteq\lambda\cdot\llg(\eta)$,
\item{(c)} If $\eta\initialeq\nu$ are both in $u$, then $M^p_\eta=M^p_\nu\rest
\lambda\cdot\llg(\nu)$,
\item{(d)} [$w${\em-closure}] $\eta\in u\,\,\&\,\,\beta\in
w[M^p_\eta]\implies\eta\rest\beta\in u$,
\item{(e)} $b^p\in [\mu]^{<\lambda}$, $c^p_\alpha\in [\lambda^+]^{<\lambda}$
for $\alpha\in b^p$,
\item{(f)${}_{\rm A}$} For $\alpha\in b^p,\iota\in c^p_\alpha$ we have
$f^\ast_\alpha\rest\iota\in u, N^p_{\alpha,\iota}\in \Gamma_\alpha$ and
$h^p_{\alpha,\iota}$
is a lawful embedding from $N^p_{\alpha,\iota}$ into
$M^p_{f^\ast_\alpha\rest\iota}$
(and hence $\card{N^p_{\alpha,\iota}}\subseteq\lambda\cdot \iota$ and $h(N^p_{\alpha,\iota})\le_K
M^p_{f^\ast_\alpha\rest\iota}$),
\item{(g)${}_{\rm A}$} If $\alpha\in b^p$ and $\iota_1<
\iota_2\in c^p_\alpha$, then
$N^p_{\alpha,\iota_1}=N^p_{\alpha,\iota_2}\rest\lambda\cdot\iota_1$ and
$h^p_{\alpha,\iota_1}=h^p_{\alpha,\iota_2}\rest N^p_{\alpha,\iota_1}$.
\end{description}
The order in $R$ is given by letting $p\le q$ iff
\begin{description}
\item{(i)} $u^p\subseteq u^q$,
\item{(ii)} for $\eta\in u^p$ we have $M^p_\eta\le_K M^q_\eta$,
\item{(iii)} $b^p\subseteq b^q$,
\item{(iv)} for $\alpha\in b^p$, we have $c^p_\alpha\subseteq c^q_\alpha$,
\item{(v)${}_{\rm A}$} for $\alpha\in b^p,\iota\in c^p_\alpha$ we have
$N^p_{\alpha,\iota}\le N^q_{\alpha,\iota}$ and $h^q_{\alpha,\iota}
\rest N^p_{\alpha,\iota}=h^p_{\alpha,\iota}$.
\end{description}

\underline{Subcase 2B}. $\lambda>\aleph_0$.
The elements of $R$ are conditions
of the form
\[
p=\langle u^p, \bar{M}^p, b^p, \bar{c}^p_\alpha,\bar{d}^p_{\alpha,\iota},
\bar{(N,h)}^{p}_{\alpha,\iota,\Upsilon}\rangle
\]
where
\begin{itemize}
\item
 $\bar{M}^p=\langle M^p_\eta:\,\eta\in u^p\rangle,$
 \item
 $\bar{c}^p_\alpha=\langle c^p_\alpha:\,\alpha\in b^p\rangle$
\item 
$\bar{d}^p_{\alpha,\iota}= \langle d^p_{\alpha,\iota}:\,
\alpha\in b^p,\iota\in c^p_\alpha\rangle$,
\item
$\bar{(N,h)}^{p}_{\alpha,\iota,\Upsilon}=
\langle (N^{p}_{\alpha,\iota,\Upsilon}, h^p_{\alpha,\iota,\Upsilon}):\,
\Upsilon\in d^p_{\alpha,\iota}, \alpha\in
b^p,\iota\in c^p_\alpha\rangle$
\end{itemize}
and
\begin{description}
\item{(a)-(e)} from Subcase 2A hold,
\item{(f)${}_{\rm B}$} for $\alpha\in b^p,\iota\in c^p_{\alpha}$ we have
$d^p_{\alpha,\iota}\in [\lambda]^{<\lambda}$,
\item{(h)} for $\alpha\in b^p,\iota\in c^p_\alpha$ we have
$f^\ast_\alpha\rest\iota\in u$ and for each $\Upsilon\in d^p_{\alpha,
\iota}$ we have
$N^p_{\alpha,\iota,\Upsilon}\in \Gamma_\alpha$ and
$h^p_{\alpha,\iota,\Upsilon}:\,N^p_{\alpha,\iota,\Upsilon}\into
M^p_{f^\ast_\alpha\rest \iota}$ is a lawful embedding (and hence
$\card{N^p_{\alpha,\iota,\Upsilon}}\subseteq\lambda\cdot i$
and $h^p_{\alpha,\iota,\Upsilon}(N^p_{\alpha,\iota,\Upsilon})\le
M^p_{f^\ast_\alpha\rest i}$),
\item{(j)} if $\alpha\in b^p$ and $\iota_1<\iota_2\in c^p_\alpha$ while
$\Upsilon\in d^p_{\alpha,\iota_2}$, then 
\[
(N^p_{\alpha,\iota_2,\Upsilon}\rest\lambda\cdot\iota_1, 
h^p_{\alpha,\iota_2,\Upsilon}\rest\lambda\cdot\iota_1)=
(N^p_{\alpha,\iota_1,\Upsilon'}, 
h^p_{\alpha,\iota_1,\Upsilon'})
\]
for some $\Upsilon'\in d^p_{\alpha,\iota_1}$.
\end{description}
The order in $R$ is given by letting $p\le q$ iff (i)-(iv) from
Subcase 2A hold and
\begin{description}
\item{(v)${}_{\rm B}$} for $\alpha\in b^p,\iota\in c^p_\alpha,\Upsilon\in
d^p_{\alpha,\iota}$ we have $N^p_{\alpha,\iota,\Upsilon}\le
N^q_{\alpha,\iota,\Upsilon'}$ and
$h^p_{\alpha,\iota,\Upsilon}\subseteq h^q_{\alpha,\iota,\Upsilon'}$
for some $\Upsilon'\in d^q_{\alpha,\iota}$.
\end{description}

If $G$ is $R$-generic, then we let for $\eta\in\TT$
\[
\Delta^{\beta,j}_\eta=\{M^p_\eta:\,p\in G\,\,\&\,\,\eta\in u^p\}.
\]

\underline{Case 3}. If $\lambda=\aleph_0$, this case does not 
occur. If $\lambda>\aleph_0$, then we are given $\alpha<\mu$ and
a $P_\beta$-name $\name{K}=\name{K}_{\ap, j}$ of a workable
$\lambda$-approximation family such that for some $j'<j$ we have had
$\name{K}_{\ap,j'}=\name{K}_{\ap,j}$ and the forcing $\name{R}^\beta_{j'}$ was
defined by Case 2. In $V^{P_\beta\ast\name{Q}^\beta_{j'}\ast
\name{R}^\beta_{j'}}$, let $G$ be the generic of $R^\beta_{j'}$ over
$V^{P_\beta\ast \name{Q}^\beta_{j'}}$
and let
$R^\beta_j\deq$
\[
\{(N,h):\,(\exists p\in G)
(\exists \iota\in c^p_\alpha)(\exists\Upsilon\in d^p_{\alpha,\iota})
[(N,h)=(N^p_{\alpha,\iota,\Upsilon}, h^p_{\alpha,\iota,\Upsilon})]\}
\]
ordered by $(N_1,h_1)\le (N_2,h_2)$ iff $N_1\le N_2$ and $h_1=h_2\rest N_1$.

\underline{Case 4}. For some $P_\beta\ast\name{Q}^\beta_j$-names
of a workable $\lambda$-approximation family $\name{K}=\name{K}^\beta_{\ap, j}$
and
a member $\name{\Gamma}^-=\name{\Gamma}^-_{\beta,j}$ of
$K_{\md}^-[\name{K}_{\ap}]$ such that 
\[
\forces_{P_{\beta}\ast\name{Q}^\beta_j}
``\{\card{M}:\,M\in \name{\Gamma}^-
\}\subseteq[\Ev]^{<\lambda}"
\]
we have (working in $V^{P_\beta}\ast\name{Q}^\beta_j$),
\[
R=\left\{
\langle M,N\rangle:\,
M,N\in K
\,\,\&\,\,M=N\rest\Ev\,\,
\&\,\,M\in \Gamma^-\right\}
\]
ordered by
\[
\langle M_1, N_1\rangle \le \langle M_2, N_2\rangle\mbox{ iff } [M_1
\le M_2\mbox{ and }
N_1
\le N_2].
\]
\end{description}
\end{Definition}

\begin{Discussion}
We now prove a series of Claims which taken together imply the
Theorem. These Claims are formulated for $\beta<\lambda^{++}, j<\mu$ and are
proved by induction on $\beta$ and $j$. Let us fix $\beta<\lambda^{++}$ and
$j<\mu$ and assume that we have arrived at the induction step for $(\beta,j)$.
We work in $V^{P_\beta\ast\name{Q}^\beta_j}$ and let $R=R^\beta_j$.
\end{Discussion}

\begin{Claim}\label{facts}
Suppose $R$ is defined by Case 2 of Definition
\ref{candidate}.

{\noindent (1)} If
$p\in R$, then
for any $\eta\in u^p$ and $\CC\subseteq u^p$ a chain with $\bigcup\CC=\eta$, we
have $M^p_\eta=\bigcup\{M^p_\nu:\,\nu\in\CC\}$.

{\noindent (2)} For every
$\eta\in \TT$ and $v\in [\lambda^+]^{<\lambda}$
such that $v\subseteq \lambda\cdot\lg(\eta)$, the set
\[
\JJ_{\eta,v}\deq\{p\in R:\,\eta\in u^p
\,\,\&\,\,v\subseteq\card{M^p_\eta}\}
\]
is a dense open subset of $R$.

{\noindent (3)}  Suppose that $p\in R$ is given and
that for some $\eta\in u^p$ we are given $M\ge M^p_\eta$ with 
$\card{M}\subseteq \lambda\cdot\llg(\eta)$.
Then there is $q \ge p$ with $M^q_\eta=M$.

{\noindent (4)} Suppose that $\lambda>\aleph_0$ and that $p\in R$,
$\alpha<\mu$
and $(N_1, h_1), (N_2, h_2)$ are such that for some $\iota_1,\iota_2\in
c^p_\alpha$ and
$\Upsilon_l\in d^p_{\alpha,\iota_l}$ we have 
\[
(N_l,h_l)=
(N^p_{\alpha,\iota_l,\Upsilon_l},h^p_{\alpha,\iota_l,\Upsilon_l})
\mbox{ for $l\in\{1,2\}$},
\]
while
\[
h_1
\rest(N_1\cap N_2)=h_2\rest(N_1\cap N_2).
\]
Then
\[
\DD=\{q\ge p:\,(\exists \iota\in c^q_\alpha)(\exists\Upsilon\in
d^q_{\alpha,\iota})\bigwedge_{l\in\{1,2\}} N^q_{\alpha,\iota,\Upsilon}\ge
N^p_{\alpha,\iota_l,\Upsilon_l}, h^q_{\alpha,\iota,\Upsilon}\supseteq
h^p_{\alpha,\iota_l,\Upsilon_l} \}
\]
is dense above $p$.

{\noindent (5)} Suppose that $\lambda>\aleph_0$ and $\alpha<\mu$.
Then for every $N\in \Gamma_\alpha$,
the set of all $p\in R$
such that for some $\iota,\Upsilon$ we have $N^p_{\alpha,\iota,\Upsilon}
\ge N$ is dense.
\end{Claim}

\begin{Proof of the Claim} (1) Obvious.

{\noindent (2)}
Clearly $\JJ_{\eta,v}$ is open, we shall show
that it is dense. Given $p\in
R$,
we shall define $q\in \JJ_{\eta,v}$ with $q\ge p$.
We do the definition in several steps.

\underline{Step I}. Let $u\deq u^p\cup\{\eta\}\cup\{\eta\cap\nu:\,\nu\in u^p\}$.

For $\sigma\in u$ we define $M_\sigma$ as follows. If for some
$\nu\in u^p$ we have $\sigma\initialeq\nu$, then let
$M_\sigma\deq M_\nu\rest(\lambda\cdot\lg(\sigma))$. Once this has been done, 
we have
$\sigma=\eta\notin u^p$ and we let
\[
M_\eta\deq\bigcup\{M_\tau:\,\tau\in u\,\,\&\,\,\tau\initial\eta\}.
\]
It has to be checked that this definition is valid, in particular that
$M_\sigma$ is well defined for $\sigma$ for which there are $\nu_1\neq\nu_2$ both in $u^p$
with $\sigma\initialeq \nu_1\cap\nu_2$. As $u^p$ is closed under
intersections, we have in this case that $\nu\deq\nu_1\cap\nu_2$ is in $u^p$ and
\[
M^p_{\nu_1}\rest(\lambda\cdot\lg(\nu))=M_\nu^p=
M^p_{\nu_2}\rest(\lambda\cdot\lg(\nu)),
\]
and hence $M_\sigma$ is well defined. Observe also that $u$ is closed under
intersections and that 
\[
\sigma\initialeq\tau\in u\implies M_\sigma=
M_\tau\rest(\lambda\cdot\lg(\sigma)),
\]
while clearly $\eta\in u$. Also note that if $\eta\in u^p$ we have $u=u^p$
and $M_\sigma=M^p_\sigma$ for $\sigma\in u$.

\underline{Step II}. For $\sigma\in u$, we find $M'_\sigma$ with
$M'_\sigma\ge_{K_{\ap,\alpha}} M_\sigma$ and $v\cap\lambda\cdot\lg(\sigma)
\subseteq \card{M'_\sigma}$,
while $\card{M'_\sigma}\subseteq\lambda\cdot\lg(\sigma)$
and such that
for $\sigma\initial \tau \in u$ we have $M'_\sigma=M'_\tau\rest \lambda\cdot\lg(\sigma)$.
This is done by
induction on $\lg(\sigma)$. Coming to $\sigma$,
if $\sigma=\bigcup\{\tau\in u:\,\tau\initial\sigma\}$, let
$M'_\sigma\deq\bigcup\{ M'_\tau: \tau\in u\,\,\&\,\,\tau\initial\sigma\}$.
Suppose otherwise and let 
\[
\delta\deq\bigcup\{\lambda\cdot\llg(\tau):\,\tau\in
u\,\,\&\,\,\tau\initial\sigma\}\mbox{ and }
M'\deq\bigcup\{M'_\tau:\,\tau\in u\,\,\&\,\,\tau
\initial\sigma\}.
\]
Since $M_\sigma\rest\delta\le M'$ by the inductive assumptions,
by the axiom of
end amalgamation in $K_{\ap}$ we have that $M'$ and $M_\sigma$ are compatible and have
a common upper bound $M_\sigma^{''}$ with the property $
\lambda\cdot\lg(\sigma)\supseteq\card{M^{''}_\sigma}
$ and such that for $\tau\in u$ with
$\tau\initial\sigma$, we have that $M^{''}_\sigma\rest\lambda\cdot\lg(\sigma)=
M'_\tau$. If $v\cap\lambda\cdot\llg(\sigma)\subseteq \card{M^{''}_\sigma}$
let $M^{'}_\sigma=M^{''}_\sigma$. Otherwise, let $f_\sigma$ be
a lawful embedding of $M^{''}_\sigma$ into $M^{'''}_\sigma$ with
$f_\sigma\rest\delta={\rm id}$ and $(M^{'''}_\sigma\setminus\delta)
\cap (M^{''}_\sigma\setminus\delta)=\emptyset$, while $\card{M^{'''}_\sigma}
\supseteq(v\cap\lambda\cdot\llg(\sigma))\setminus M^{''}_\sigma$.
Applying amalgamation to $M^{''}_\sigma, M^{'''}_\sigma$ we can find
$M^{'}_\sigma$ as required.

\underline{Step III}. Now for $\sigma\in u$ let $M^q_\sigma\deq M'_\sigma$.
Let
\[
u^q\deq u\cup\{\sigma\rest\beta:\,\sigma\in u\,\,\&\,\,\beta\in w[M^q_\sigma]\}.
\]
For $\sigma\in u^q$, let $M^q_\sigma\deq M'_\sigma$ if $\sigma\in u$,
and otherwise let $M^q_\sigma=M^q_\tau\rest(\lambda\cdot\lg(\sigma))$ for
any $\tau\in u$ with $\sigma\initialeq\tau$.

\underline{Subcase A}. $\lambda=\aleph_0$. Let
\[
q\deq\langle u^q,\langle M^q_\sigma:\,\sigma\in u^q\rangle, b^p,
\langle c^p_\alpha:\,\alpha\in b^p\rangle,
\langle(N^p_{\alpha,\iota}, h^p_{\alpha,\iota}):\,\alpha\in b^p,\iota\in
c^p_\alpha\rangle
\rangle.
\]

\underline{Subcase B}. $\lambda>\aleph_0$
Let
\[
q\deq\langle u^q,\langle M^q_\sigma:\,\sigma\in u^q\rangle, b^p,
\langle c^p_\alpha:\,\alpha\in b^p\rangle,
\langle d^p_{\alpha,\iota}:\,\alpha\in b^p,\iota\in c^p_\alpha\rangle,
\bar{(N,h)}^p\rangle.
\]
It is easily seen that $q$ is as required.

{\noindent (3)} All coordinates
of $q$ will be the same as the corresponding coordinates of $p$,
with the possible exception of $u^p$ and $\langle M^p_\sigma:\,\sigma\in
u^p\rangle$. Let $u=u^p$. For $\sigma\in u$ we define $M^q_\sigma\ge
M^p_\sigma$ by induction on $\llg(\sigma)$. The inductive hypothesis
is that if $\sigma\initialeq\eta$, then
$M^q_\sigma=M\rest\lambda\cdot\llg(\sigma)$.
The proof is similar to that of of Step II of (2).
Coming to $\sigma$, let $\delta$ and $M'$ be defined as there. If
$\sigma\initialeq\eta$, then we let $M^q_\sigma
=M\rest\lambda\cdot\llg(\sigma)$,
and we have $M^q_\sigma\ge M^p_\sigma$ by the choice of $M$.

Otherwise, let $\nu=\sigma\cap\eta$, and hence $\nu\in u$ and $\llg(\nu)
<\llg(\sigma)$. We have that 
\[
M^q_\nu=M\rest\lambda\cdot\llg(\eta)\ge
M^p_\nu=M^p_\sigma\rest\lambda\cdot\llg(\sigma).
\]
Hence we can find $M^q_\sigma\ge M^p_\sigma$ with
$M^q_\sigma\rest\lambda\cdot\llg(\nu) =M^q_\nu$ by end amalgamation.

Once the induction is done, we define $u^q$ exactly as in the Step III
of (2), and define $M^q_\sigma$ for $\sigma\in u^q$ accordingly.

{\noindent (4)}  Let $r\ge p$, we shall find
$q\in \DD$ with $q\ge r$. For some $\iota'_1,\iota'_2,\Upsilon'_1,\Upsilon'_2$
we have $N_1\le N^r_{\alpha,\iota'_1,\Upsilon'_1}$ and
$N_2\le N^r_{\alpha,\iota'_2,\Upsilon'_2}$, and similarly for $h_1,h_2$. For
simplicity of notation,
we assume $\iota'_l=\iota_l$ and $\Upsilon'_l=\Upsilon_l$ for $l\in\{1,2\}$.
Let $\iota=\max\{\iota_1,\iota_2\}$ and
let $\eta=f^\ast_\alpha\rest\iota$. Let $M'\deq M^r_\eta$,
which is well defined. Hence
$h_1$ and $h_2$ are lawful embeddings of $N_1$ and $N_2$ into
$M'$ respectively.

First define a lawful isomorphism $g_0$ from $M'$ onto some $M_1$ such that for
$x \in N_l$ we have $g(h_l(x))=x$ for $l\in\{1,2\}$. This is possible
because $h_1$ and $h_2$ agree on $N_1\cap N_2$. Hence we have that $N_1\in
\Gamma$ and $N_1\le M_1$. As $\Gamma\in K_{\rm md}$, there is a lawful
embedding $g_1:\,M_1\into M_2$ for some $M_2\in \Gamma$ such that
$g_1$ is the identity on $N_1$. Without loss of generality, again as $\Gamma\in
K_{\rm md}$, we can assume that $\card{M_2}\subseteq\lambda\iota$.

Now let $g_2$ be a lawful isomorphism between $M_2$ and $M_3$ such that $g_2
\rest N_1={\rm id}$ and $g_2(g_1(x))=x$
for $x\in N_2$. Then $N_2\le M_3$, so we can find $M_4\in \Gamma$ and
a lawful embedding $g_3:\,M_3\into M_4$ over $M_3$. Then $N_1, N_2\le M_4$.
Without loss of generality we have $\card{M_4}\subseteq\lambda\cdot\iota$.
Finally, there is a lawful isomorphism $g$ between $M_4$ and
some $M$ such that $g\rest N_l=h_l$ for $l\in\{1,2\}$ and $M'\le M$.
By (3) we can find $q'\ge r$ such that $M^{q'}_\eta=M$.
We shall define $q\ge q'$ so that all the coordinates of $q$ are the same as
the corresponding coordinates of $q'$, except that we in addition choose some
\[
\Upsilon\in
\lambda\setminus (d^{q'}_{\alpha,\iota}\cup\bigcup\{d^q_{\alpha,j}:\,j\in
w[M_4]\}
\]
and let $N^q_{\alpha,\iota,\Upsilon}= M_4$, while $h^q_{\alpha,\iota,\Upsilon}=g$.
We let $N^q_{\alpha,j,\Upsilon}=M_4\rest\lambda\cdot j$ for $j\in w[M_4]$, and
similarly for $h^q_{\alpha,j,\Upsilon}$.
Then $q$ is as required.

{\noindent (5)} Similar to (4), using (3) and (4).
$\eop_{\ref{facts}}$
\end{Proof of the Claim}

\begin{Claim}\label{A} $(1)$ If $\lambda>\aleph_0$, then
$R$ is a $<\lambda$-complete forcing.

{\noindent (2)}
If $R$ was defined by one of the Cases 3-4 or by Case 2 and $\lambda=\aleph_0$,
then
every increasing sequence
in $R$ of length $<\lambda$
has a least upper bound.

{\noindent $(3)$} Suppose $R$ was defined by Case 2.
Then,
for every $\alpha<\mu$
\[
\bigcup_{\gamma<\lambda^+}\Delta^{\beta,j}_{f^\ast_\alpha\rest \gamma}\in
K_{\md}^-[K_{\ap,\alpha}]
\]
holds in $V^{P_\beta\ast\name{Q}^\beta_j\ast\name{R}^\beta_j}$. In addition,
$\bigcup\{\card{M}:\,M\in \bigcup_{\gamma<\lambda^+}
\Delta^{\beta,j}_{f^\ast_\alpha\rest \gamma}
\}=\lambda^+$.  

If $\langle \eta_i:\,i<i^\ast<\lambda
\rangle$ is a $\unlhd$-increasing sequence of elements of $\TT$,
and $M_i\in \Delta^{\beta,j}_{\eta_i}$ for $i<i^\ast$,
then $\bigcup_{i<i^\ast} M_i\in
\Delta^{\beta,j}_{\bigcup_{i
<i^\ast}\eta_i}$.

\end{Claim}

\begin{Proof of the Claim} 
$(1)$ If $R$ is defined by Case 1, this follows by the definition of that
case. For Cases 3-4 the conclusion follows by (2). We give the proof for
Case 2. We deal with the situation $\lambda>\aleph_0$. The other case is
trivial and in that case we actually obtain the existence of lubs.

Suppose that $\bar{q}=\langle q_i:\,i<i^\ast<\lambda\rangle$ is an increasing
sequence in $R$. Without loss of generality, $i^\ast$ is a limit ordinal.
Let $b\deq\bigcup_{i<i^\ast} b^{q_i}$ and for $\alpha\in b$ let
$c_\alpha=\bigcup\{c^{q_i}_\alpha:\,i<i^\ast\,\,\&\,\,\alpha\in b^{q_i}\}$. Let
$u=\bigcup_{i<i^\ast}u^{q_i}$ and for every $\eta\in u$ let
$M_\eta=\bigcup_{i<i^\ast,\eta\in u^{q_i}}M^{q_i}_\eta$.

Let $\theta=\card{i^\ast}$, and so $\theta<\lambda$. For $\alpha\in b$ and
$\iota\in c_\alpha$, $\Upsilon\in \bigcup_{i<i^\ast} d^{q_i}_{\alpha,\iota}$
we let
\[
(N_{\alpha,\iota,\theta\cdot i+\Upsilon}, h_{\alpha,\iota,\theta\cdot
i+\Upsilon})= (N^{q_i}_{\alpha,\iota,\Upsilon}, h^{q_i}_{\alpha,\iota,\Upsilon})
\]
if this is defined. Let 
\[
d_{\alpha,\iota}=\{\theta\cdot i+\Upsilon:\,
i<i^\ast\mbox{ and } N_{\alpha,\iota,\theta\cdot i+\Upsilon}\mbox{ well
defined}\}.
\]
Let
\[q=\langle u, \langle M_\eta:\,\eta\in u\rangle, b,
\langle c_\alpha:\,\alpha\in b\rangle,\langle d_{\alpha,\iota}:\,
\alpha\in b,\iota\in c_\alpha\rangle,\bar{(N,h)}^q\rangle
\]
where $\bar{(N,h)}^q=
\langle (N_{\alpha,\iota,\Upsilon},h_{\alpha,\iota,\Upsilon}):\,\Upsilon\in
d_{\alpha,\iota},\iota\in c_\alpha,\alpha\in b\rangle$.
It is easily seen that $q$ is an upper bound of all $q_i$ (although not a
least upper bound, which may not exist).

{\noindent (2)} The Case 2, Subcase $\lambda=\aleph_0$ is trivial.
We distinguish Cases 3 and 4, according to Definition \ref{candidate}.

Suppose that $\bar{q}=\langle q_i:\,i<i^\ast<\lambda\rangle$
is an increasing sequence in $R$. We shall
define the lub $q$ of $\bar{q}$.

\underline{Case 3}. Let $j'<j$, $\alpha<\mu$ and $G$ be as in the definition
of the forcing and let $q_i=(N_i,h_i)$ and $p_i\in G$ for $i<i^\ast$ be such
that
$\alpha\in b^p$ and for some $\iota_i\in c^p_\alpha$ and $\Upsilon_i\in
d^p_{\alpha,\iota}$ we have
\[
(N_i,h_i)=(N^{p_i}_{\alpha,\iota_i,\Upsilon_i},
h^{p_i}_{\alpha,\iota_i,\Upsilon_i}).
\]
Without loss of generality, $i^\ast$ is a limit ordinal. As we know by (1) of
the induction hypothesis that $R^\beta_{j'}$ is $(<\lambda)$-complete, there is
$p\in G$ with $p\ge p_i$ for all $i<i^\ast$. In
$V^{P_\beta\ast\name{Q}^\beta_{j'}}$ let $\eta=\bigcup\{\nu\in u^p:\,\nu\initialeq
f^\ast_\alpha\}$ and let $M=\bigcup\{M_\nu:\,\nu\in u^p\,\,\&\,\,\nu\initialeq
f^\ast_\alpha\}$. Letting $N=\bigcup_{i<i^\ast} N_i$ and
$h=\bigcup_{i<i^\ast} h_i$ we obtain that $h$ is a lawful embedding from $N$ into
$M$. As $\Gamma_\alpha\in K_{\rm md}$, by the $(<\lambda)$-closure of
$\Gamma_\alpha$ we obtain $N\in \Gamma_\alpha$. Consider the set $\DD$ defined
by
\[
\{q\ge p:\,\eta\in u^q\,\,\&\,\,M_\eta^q\ge M \,\,\&\,\,
(\exists \iota\in c_\alpha^q)(\exists\Upsilon\in
d^q_{\alpha,\iota})N^q_{\alpha,\iota,\Upsilon}=N\,\,\&\,\,h^q_{\alpha,\iota}=h\}.
\]
Note that $\DD\in V^{P_\beta\ast\name{Q}^\beta_{j'}}$ as all forcings involved
are
$(<\lambda)$-closed by the induction hypothesis. 

\begin{Subclaim}\label{densityD} $\DD$ is dense above $p$ in the forcing
$R^\beta_{j'}$.
\end{Subclaim}

\begin{Proof of the Subclaim} Let $r\ge p$ be given. For every $\nu\in u^p$
with $\nu\initialeq f^\ast_\alpha$ we have $M^r_\nu\ge M^p_\nu$.
Hence if $\eta\in u^r$ we have $M^r_\eta=\bigcup\{M^r_\nu:\,\nu\initialeq\eta\}\ge
M$. By Claim \ref{facts}(2) we can without loss of generality assume
that this is the case. Let $u^q=u$ and for $\nu\in u^q$ let $M^q_\nu=
M^r_\nu$.

Let $\iota=\llg(\eta)$, and hence
$\card{M}\subseteq\lambda\cdot\iota$. Further let
$c^q_\alpha=c^r_\alpha\cup\{\iota\}$. Let $\Upsilon$ be such that
$\Upsilon\notin\bigcup_{\iota'\le\iota} d^r_{\alpha,\iota'}$ and for $\iota'\le
\iota$ with $\iota'\in c^q_\alpha$ let 
\[
(N^q_{\alpha,\iota',\Upsilon}, 
h^q_{\alpha,\iota',\Upsilon})=(N\rest\lambda\cdot\iota', h\rest\lambda\cdot\iota').
\]
We complete the definition of $q$ in the obvious fashion.
Hence $q\ge r$ and $q\in \DD$.
$\eop_{\ref{densityD}}$
\end{Proof of the Subclaim}

By the Subclaim it follows that there is $q\in G\cap\DD$, and this $q$
witnesses that $(N,h)\in R$. Obviously, $(N,h)$ is the lub
of $\langle (N_i,h_i):\,i<i^\ast\rangle$.

\underline{Case 4}. If $\langle \langle M_i,
N_i\rangle:\,i<i^\ast<\lambda\rangle$ is
increasing in $R$ then clearly $\langle \bigcup_{i<i^\ast} M_i,\cup_{i<i^\ast}
N_i
\rangle$ is the lub.

{\noindent (3)} 
We first prove the second statement.
So, let $\langle \eta_i:\,i<i^\ast\rangle$ be as in the Claim.
Without loss of generality, $i^\ast$ is a limit ordinal. Let
$\eta\deq\bigcup_{i<i^\ast}\eta_i$. For $i<i^\ast$ let $p_i\in G=G_{R^\beta_j}$
be such that $\eta_i\in u^{p_i}$ and $M_i=M^{p_i}_{\eta_i}$.
Let $p$ be an upper bound of $\langle p_i:\,i<i^\ast\rangle$ with
$p\in G$, which exists by (1).
Now let $q\in G$ be such that
$\eta\in u^q$ and $p\le q$, which exists by Claim \ref{facts}. Note that
we have $M^q_\eta=\bigcup_{i<i^\ast}
M^q_{\eta_i}$. Let now $r$ be defined by $u^r=u^p\cup\{\eta\}$, and for
$\nu\in u^p$ we have $M^r_\nu=M^p_\nu$, while $M^r_\eta=\bigcup_{i<i^\ast}
M_i$. We also redefine $c_\alpha, d_{\alpha,\iota}$
and $\bar{(N,h)}^q$ to accommodate the fact that we have shrunk
$u^q$, for example by using the corresponding coordinates of $p$.
This gives us a well defined condition $r$. We now claim that
$r\le q$. We only need to show that $M^r_\eta\le M^q_\eta$, which follows
by Definition \ref{general} (1)(c). As $G$ is generic, and $q\in
G$,
we have $r\in G$.

For the first statement, suppose that $M_i\,(i<i^\ast<\lambda)$ are in
$\bigcup_{\gamma<\lambda^+}\Delta^{\beta,j}_{f^\ast_\alpha\rest \gamma}$ and
let $\gamma_i$ for $i<i^\ast$
be such that $M_i\in \Delta^{\beta,j}_{f^\ast_\alpha\rest \gamma}$.
Let $\eta_i\deq f^\ast_\alpha\rest \gamma_i$ for $i<i^\ast$.
Let $\eta\deq\bigcup_{i<i^\ast}\eta_i$. Now proceed as above.
This proves that $\bigcup_{\gamma<\lambda^+}\Delta^{\beta,j}_{f^\ast_\alpha\rest \gamma}$ is
a $(<\lambda)$-closed subset of $K_{\ap}$ and it is equally easy
to see that it is directed. To see that
$\bigcup\{\card{M}:\,M\in
\bigcup_{i<\lambda^+}\Delta^{\beta,j}_{f^\ast_\alpha\rest \gamma}\}=\lambda^+$, apply Claim
\ref{facts}, and this of course implies that (iii) of Definition
\ref{md}(2) holds.

$\eop_{\ref{A}}$
\end{Proof of the Claim}

\begin{Notation} The upper bound $q$ of $\bar{q}$ that is constructed as in
the proof of Claim \ref{A}(1) will be called a canonical upper bound (cnub)
of $\bar{q}$.
\end{Notation}

\begin{Note}\label{easyr} The same
proof given above shows that if  $\name{R}^\beta_{j}$ is
defined by Case 2 of Definition \ref{candidate}
and $\eta\in \TT$, then in
$V^{P_{\beta}\ast\name{Q}^\beta_j\ast\name{R}^\beta_j}$ we have that
$\Delta_\eta^{\beta,j}$ is an element of 
$K_{\md}^-[K_{\ap}^{\lambda\cdot\lg(\eta)}]$, where
$K_{\ap}=K_{\ap}^{\beta,j}$. \end{Note}

\begin{Claim}\label{astepsilon}
If $\lambda>\aleph_0$, then
$R$ satisfies
$\ast^\varepsilon_\lambda.$

If $\lambda=\aleph_0$, then $R$ satisfies ccc.
\end{Claim}

\begin{Proof of the Claim} We distinguish various cases of Definition
\ref{candidate}.

\underline{Case 1}. $R$ is defined by Case 1 of Definition \ref{candidate}.
The conclusion follows by the assumptions.

\underline{Case 2}. (main case) $R=R^\beta_j$ is defined by Case 2 of
Definition \ref{candidate}. As Subcase B is more difficult, we start by it.

\underline{Subcase B}. $\lambda>\aleph_0$. Let $K=K^{\beta,j}_{\ap}$, and
let us follow the rest of the notation of Definition \ref{candidate} as well.
By our assumptions we have $\card{\TT}=\lambda^+$ and by
Claim \ref{A}, the equality $(\lambda^+)^{<\lambda}=
\lambda^+$ holds. Also, for every $j<\lambda^+$
we have that $K\rest j\deq
\{M\in K:\,\card{M}\subseteq j\}$ has cardinality $\le\lambda$.

We first define several auxiliary functions. Let $g^\ast_0:\,\TT\into
\lambda^+$ be a bijection and let
\[
g^\ast_1:\,\lambda^+\into{}^\lambda 2
\]
be a 1-1 function.

\begin{Subclaim}\label{stublic} There is a function 
$g^\ast_2:\,K\into\lambda$
such that for every $N_1, N_2\in K$ we have
\[
g_2^\ast(N_1)=g_2^\ast(N_2)\,\,\&\,\,\beta\in w[N_1]\cap w[N_2]
\implies N_1\rest\lambda\beta= N_2\rest\lambda\beta.
\]
\end{Subclaim}

\begin{Proof of the Subclaim}
For $N\in K$ define
\[
o[N]\deq\{\lambda\gamma:\,\gamma\in w[N]\}\cup\card{N},
\]
\[
\xi[N]\deq\min\{\zeta<\lambda:\,\beta\neq\gamma\in o[N]\implies
g^\ast_1(\beta)\rest\zeta\neq g^\ast_1(\gamma)\rest\zeta\},
\]
\[
\Xi_N\deq\{g^\ast_1(\beta)\rest \xi[N]:\,\beta\in o[N]\}.
\]
Note that $\xi[N]$ is well defined because $g^\ast_1$ is 1-1 and $\card{o[N]}
<\lambda$.

For $\alpha<\lambda^+$, let $g_\alpha:\,\alpha\into\lambda$ be one to one. 
For $N\in K$
let ${\frak A}_N$ be a model with universe included in $\Xi_N$ such that the
function
\[
\beta\mapsto g^\ast_1(\beta)\rest \xi[N]
\]
is an isomorphism from $N$ onto ${\frak A}_N$. Let $<_N$ be a well ordering of
$\Xi_N$ such that $\beta\mapsto g^\ast_1(\beta)\rest \xi[N]$ is
an isomorphism from $(o[N], <)$ onto $(\Xi_N, <_N)$. 
Let 
\[
R^N\deq\{(g^\ast_1(\beta)\rest\xi[N],g^\ast_1(\gamma)\rest\xi[N],
g_\gamma(\beta)):\,\beta<\gamma\mbox{ both in }o[N]\}.
\]
Notice that $(i,j,k_1),(i,j,k_2)\in R^N\implies
k_1=k_2$ by the choice of $\xi[N]$.
Let $\HH^\ast:\,\HH(\lambda)\into\lambda$ be one to one, which exists as
$\lambda^{<\lambda}= \lambda$. We define 
\[
g^\ast_2(N)\deq\HH^\ast(\langle \xi[N], \Xi_N,
{\frak A}_N, <_N, R^N\rangle).
\]
Clearly $g^\ast_2$ is a well defined function from $K$ to $\lambda$.
Let us show that it has the required properties.

Suppose $g^\ast_2(N_1)=g^\ast_2(N_2)$ and $\beta^\ast\in w[N_1]\cap w[N_2]$.
Firstly, we have that $\xi[N_1]=\xi[N_2]=\xi$ and the functions
\[
f_1:\,\beta\mapsto g^\ast_1(\beta)\rest \xi[N_1]\mbox{ for }\beta\in o[N_1]
\]
and
\[
f_2:\,\beta\mapsto g^\ast_1(\beta)\rest \xi[N_2]\mbox{ for }\beta\in o[N_2]
\]
are one to one and onto the same set $\Xi_{N_1}=\Xi_{N_2}=\Xi$. Furthermore,
both $f_1$ and $f_2$ are
order preserving and $<_{N_1}=<_{N_2}$. Hence there is a
one to one $<$-preserving function $g:\,o[N_1]\into o[N_2]$
given by $g(\beta)=f_2^{-1}(f_1(\beta))$.

We claim that for every 
$\beta\in w[N_1]\cap w[N_2]$ we have
$g(\lambda\beta)=\lambda\beta$. Namely suppose not,
say $g(\lambda\beta)=\gamma$ and $\lambda\beta<\gamma$. Then
$f_2(g(\lambda\beta))=f_2(\gamma)>f_2(\lambda\beta)$, and hence
$f_1(\lambda\beta)>f_2(\lambda\beta)$, which means that $g^\ast_1(\lambda\beta)
\rest\xi>_{N_1}g^\ast_1(\lambda\beta)\rest\xi$, a contradiction. A similar
contradiction can be obtained by assuming that $g(\lambda\beta)<\lambda\beta$.

If $\gamma\in N_1\rest\lambda\beta^\ast$ then $g(\gamma)<g(\lambda\beta^\ast)=
\lambda\beta^\ast$. By the definition of $g$ we have 
$g^\ast_1(\gamma)\rest\xi =g^\ast_1(g(\gamma))\rest\xi$. Hence,
$\lambda\beta^\ast,g(\gamma)\in o[N_2]$ and $g(\gamma)<\lambda\beta^\ast$.
So $(g_1^\ast(g(\gamma))\rest\xi, g^\ast_1(\lambda\beta^\ast)\rest\xi,
g_{\lambda\beta^\ast}(g(\gamma)))\in
R^{N_2}=R^{N_1}$. As also
\[
(g_1^\ast(\gamma)\rest\xi, g^\ast_1(\lambda\beta^\ast)\rest\xi,
g_{\lambda\beta^\ast}(\gamma))\in
R^{N_1},
\]
we have that $g_{\lambda\beta^\ast}(g(\gamma))=
g_{\lambda\beta^\ast}(\gamma)$ and hence $g(\gamma)=\gamma$. In particular
$\gamma\in o[N_2]$. As we have $\gamma\in N_1$, we have
$g^\ast_1(\gamma)\rest\xi\in {\frak A}_{N_1}$, and hence
$g^\ast_1(\delta)\rest\xi=g^\ast_1(\gamma)\rest\xi=g^\ast_1(g(\gamma))\rest\xi$
for some $\delta\in N_2$. 
As $\xi=\xi[N_2]$, we have that $\delta=g(\gamma)=\gamma$, so $\gamma\in
N_2\rest\lambda\beta^\ast$.

This argument shows that $N_1\rest\lambda\beta^\ast\subseteq
N_2\rest\lambda\beta^\ast$, and it can be shown similarly that
$N_1\rest\lambda\beta^\ast=N_2\rest\lambda\beta^\ast$ as sets
and as models.
$\eop_{\ref{stublic}}$
\end{Proof of the Subclaim}

For $p\in R$, let $\langle (\eta(p,i):\,i<i(p)\rangle$ list $u^p$
with no repetitions, and let $\xi(p)$ be the minimal $\xi<\lambda$
such that
\[
\langle g_1^\ast(g^\ast_0(\eta(p,i)))\rest\xi:\,i<i(p)\rangle
\]
is without repetitions (which exists as $g^\ast_0$ and $g^\ast_1$ are 1-1
and $\langle (\eta(p,i):\,i<i(p)\rangle$ is without repetitions). Let
\[
g^\ast_3:\,R\into\lambda
\]
be such that for $p,q\in R$ with $g^\ast_3(p)=g^\ast_3(q)$
we have
\begin{description}
\item{(a)} $i(p)=i(q)$,
\item{(b)} the mapping defined by sending $\eta(p,i)\mapsto \eta (q,i)$
preserves
\[
``\nu\initialeq\eta", ``\neg (\nu\initialeq\eta)", ``\nu_1\cap\nu_2=\nu",
``\neg (\nu_1\cap\nu_2=\nu)",
\]
\item{(c)} $\xi(p)=\xi(q)$,
\item{(d)} for $i<i(p)$ we have $g^\ast_1(g^\ast_0(\eta(p,i)))\rest\xi(p)=
g^\ast_1(g^\ast_0(\eta(q,i)))\rest\xi(p)$ (recall that
$\card{{}^{\lambda>}2}=\lambda$),
\item{(e)} for $i<i(p)$ we have $g_2^\ast(M^p_{\eta(p,i)})=
g_2^\ast(M^q_{\eta(q,i)})$.
\end{description}

The existence of such a function can be shown by counting.

\begin{Subclaim}\label{film} If $g^\ast_3(p)=g^\ast_3(q)$, then the mapping
sending $\eta(p,i)$
to $\eta(q,i)$ for $i<i(p)=i(q)$, is the identity on
$u^{p}\cap u^{q}$.
\end{Subclaim}

\begin{Proof of the Subclaim}
Suppose that $\eta\in u^p\cap u^q$. Let $i$ be such 
that $\eta=\eta(p,i)$. Letting $\xi\deq\xi(p)=
\xi(q)$,we have 
\[
g^\ast_1(g^\ast_0(\eta))\rest\xi=g^\ast_1(g^\ast_0(\eta(q,i)))\rest\xi.
\]
By the definition of
$\xi$ and the fact that $\eta\in u^q$, we must have $\eta(q,i)=\eta$.
$\eop_{\ref{film}}$
\end{Proof of the Subclaim}

Let us also fix a
bijection
\[
F:\,\lambda\times
{}^{\lambda>}([\lambda^+]^{<\lambda})\into\lambda^+
\]
and let $C$ be a club of $\lambda^+$ such that for every $j\in 
S^{\lambda^+}_\lambda\cap C$ we have
\[
\beta<\lambda
\,\,\&\,\,u\in{}^{\lambda>}([j]^{<\lambda})\implies
F\left((\beta,u)\right)<j.
\]
We describe a winning strategy for I in
$\ast^\varepsilon_\lambda[R]$.
Given $0<\zeta<\varepsilon$ and suppose that
\[
\left\langle(\langle q^\xi_s:\,s<\lambda^+\rangle, f_\xi),
\langle p^\xi_s:\,s<\lambda^+\rangle:\,\xi<\zeta\right\rangle
\]
have been played so far and I has played according to
the strategy. By Claim \ref{A}(1),
we can let player
I choose $q^\zeta_s$ as a cnub of $\langle p^\xi_s:\,\xi<\zeta
\rangle$.
Next we describe the choice of $f_\zeta$. Let
$C_\zeta\deq C$ and
define $g_\zeta$
which to an ordinal $j\in S^{\lambda^+}_\lambda$ assigns
\[
(g^\ast_3(q^\zeta_j), 
\langle w[M^{q^\zeta_j}_{\eta(q^\zeta_j,i)}]\cap j
:\, i<i(q^\zeta_j)\rangle).
\]
Then let
\[
f_\zeta\deq (F\circ g_\zeta)\rest (C_\zeta\cap
S^{\lambda^+}_\lambda)\cup 0_{\lambda^+\setminus (C_\zeta\cap
S^{\lambda^+}_\lambda)}.
\] 
Let us check that this definition is as required. It follows from the choice of
$C$ that
each $f_\zeta$ is regressive on $C_\zeta\cap
S^{\lambda^+}_\lambda$. Let $E\subseteq C$
be a club of $\lambda^+$ such that
\[
[j\in E\cap S^{\lambda^+}_\lambda\,\,\&\,\,j'<j]\implies
(\forall\zeta<\varepsilon)(\forall i<i(q^\zeta_{j'}))\,
[w[M^{q^\zeta_{j'}}_{\eta(q^\zeta_{j'},i)}]\subseteq j].
\]
Suppose that $j'<j\in E\cap S^{\lambda^+}_\lambda$ are such that
\[
\bigwedge_{\zeta<\varepsilon}f_\zeta(j')=f_\zeta(j).
\]
We define an upper bound to
\[
\{p^\zeta_{j'}:\,\zeta<\varepsilon\}\cup\{p^\zeta_j:\,\zeta<\varepsilon\}.
\]
As we have $q^{\zeta+1}_s\ge p^\zeta_s$ for all
$\zeta<\varepsilon$ and $s<
\lambda^+$, and $\varepsilon$ is a limit ordinal, it suffices to
define an upper bound to
\[
\{q^\zeta_{j'}:\,\zeta<\varepsilon\}\cup\{q^\zeta_j:\,\zeta<\varepsilon\}.
\]
We first define $q_l$ as a cnub of $\{q^\zeta_l:\,\zeta<\varepsilon\}$
for $l\in \{j',j\}$, and
we shall now describe an upper bound $r$ of $q_{j'}$ and $q_j$. Notice
that $u^{q_l}=\bigcup_{\zeta<\varepsilon}u^{q_l^\zeta}$ for $l\in\{j',j\}$.

Let
\[
u \deq u^{q_{j'}}\cup u^{q_j}\cup\{\eta\cap\nu:\,\eta\in u^{q_{j'}}\,\,\&\,\,\nu\in
u^{q_j}\}.
\]
Clearly $\card{u}<\lambda$ and $u$ is closed under intersections.
For $\eta\in u$, let 
\[
M_\eta\deq M^{q^\zeta_l}_\nu\rest\lambda\cdot\lg(\eta)
\]
for
any $l\in \{j',j\}$, $\zeta<\varepsilon$ and $\nu\in u^{q^\zeta_l}$ for which
$\eta\initialeq\nu$. 

\begin{Subclaim}\label{ples} For $\eta\in u$
the model $M_\eta$ is well
defined and $\card{M_\eta}\subseteq\lambda\cdot\llg(\eta)$. For every
$l\in\{j,j'\}$ for which $\eta\in u^{q_l}$ we have $M^{q_l}_\eta=M_\eta$.
\end{Subclaim}

\begin{Proof of the Subclaim}
Firstly, note that for any $\eta\in u$ we have $\eta\initialeq\nu$ for some
$\nu\in \bigcup_{l\in\{j',j\},
\zeta<\varepsilon}u^{q^\zeta_l}$. Suppose $\eta\initialeq
\nu_1,\nu_2$ for some $\nu_1,\nu_2\in \bigcup_{l\in\{j',j\},
\zeta<\varepsilon}u^{q^\zeta_l}$ such that $\nu_k \in u^{q^{\zeta_k}_{l_k}}$
for $k\in\{1,2\}$, and $M^{q^{\zeta_1}_{l_1}}_{\nu_1}\rest\lambda\cdot\lg(\eta)
\neq M^{q^{\zeta_2}_{l_2}}_{\nu_2}\rest\lambda\cdot\lg(\eta)$. By taking the
larger of $\zeta_1,\zeta_2$, we may assume that $\zeta_1=\zeta_2=
\zeta$. By the closure under intersections, we can also assume that $l_1\neq
l_2$, so without loss of generality we have $l_1=j'$ and $l_2=j$. Let $\beta\le\lg(\eta)$ be 
minimal such that
\[
M^{q^{\zeta}_{j'}}_{\nu_1}\rest\lambda\cdot\beta
\neq M^{q^{\zeta}_{j}}_{\nu_2}\rest\lambda\cdot\beta.
\]
By the minimality of $\beta$, we have $\beta=\gamma+1$ for some $\gamma\in
w[M^{q^{\zeta}_{j'}}_{\nu_1}] \cap
w[M^{q^{\zeta}_{j}}_{\nu_2}]$. As $j\in E$ we have that
$w[M^{q^{\zeta}_{j'}}_{\nu_1}]\subseteq j$, so $\gamma\in
w[M^{q^{\zeta}_{j}}_{\nu_2}]\cap j$.
As $f_\zeta(j)=f_\zeta(j')$, there is $\nu\in u^{q^\zeta_{j'}}$
such that for some $i<i(q^\zeta_j)=i(q^\zeta_{j'})$ we have
$\nu_2=\eta(q^\zeta_{j},i)$, $\nu=\eta(q^\zeta_{j'},i)$ and
$w[M^{q^\zeta_j}_{\nu_2}]\cap j= w[M^{q^\zeta_{j'}}_{\nu}]\cap j'$.
Hence we have $g^\ast_2(M^{q^\zeta_{j'}}_\nu)=g^\ast_2(M^{q^\zeta_{j}}_{\nu_2})$,
and as $\gamma,\beta\in w[M^{q^\zeta_{j'}}_\nu]
\cap
w[M^{q^{\zeta}_{j}}_{\nu_2}]$, we have $M^{q^{\zeta}_{j'}}_{\nu}\rest\lambda\beta
=M^{q^{\zeta}_{j}}_{\nu_2}\rest\lambda\beta$. We have not arrived at a
contradiction yet, as we do not know the relationship between
$\nu$ and $\nu_1$.

As $\beta\le\lg(\eta)$, we have $\rho\deq\nu_1\rest\beta=\nu_2\rest\beta$.
Since $\beta\in w[M^{q^\zeta_{j'}}_{\nu_1}]$, we have $\rho=\nu_1\rest\beta\in
u^{q^\zeta_{j'}}$ and similarly $\rho\in 
u^{q^{\zeta}_j}$. Let $o$ be such that $\rho=\eta(q^{\zeta}_{j'}, o)$. 
By Subclaim \ref{stublic} we have
$\rho=\eta(q^{\zeta}_j, o)$. Since we have $\rho\le\nu_2$ by the choice of
$g^\ast_3$, we have $\rho\initialeq\nu$. 
So $\rho\initialeq\nu_1\cap\nu$, and as we have $\nu_1\cap\nu\in
u^{q^\zeta_{j'}}$, we obtain
\[
M^{q^{\zeta}_{j'}}_{\nu_1}\rest\lambda\beta=
M^{q^{\zeta}_{j'}}_{\nu_1}\rest\lambda\lg(\rho)=
M^{q^{\zeta}_{j'}}_{\nu}\rest\lambda\lg(\rho)=
M^{q^{\zeta}_{j}}_{\nu_2}\rest\lambda\lg(\rho)=
M^{q^{\zeta}_{j}}_{\nu_2}\rest\lambda\beta,
\]
a contradiction. This proves the first part of the statement. If $\eta\in u$
and $l\in \{j,j'\}$ is such that $\eta\in u^{q_l}$, then we have
$M_\eta=M^{q_l}_\eta$, as is clear from the definition.
$\eop_{\ref{ples}}$
\end{Proof of the Subclaim}

Now we let
\[
u^r\deq u\cup\{w[M_\sigma]:\,\sigma\in u\},
\]
and define $M^r_\sigma$ for $\sigma\in u^r$ accordingly, which is
done as in Step III of the Proof of Claim \ref{facts}(2).

We let $b^r=b^{q_j}\cup b^{q_{j'}}$ and for $\alpha\in b$ we let
$c^r_\alpha=c^{q_j}_\alpha\cup c^{q_{j'}}_\alpha$. For $\alpha\in b,
\iota\in c^r_\alpha$ we let
$d^r_{\alpha,\iota}=\{2\Upsilon:\,\Upsilon\in d^{q_{j'}}_{\alpha,\iota}\}
\cup \{2\Upsilon+1:\,\Upsilon\in d^{q_j}_{\alpha,\iota}\}$  and 
\[
(N^r_{\alpha,\iota,2\Upsilon}, h^r_{\alpha,\iota,2\Upsilon})=
(N^{q_{j'}}_{\alpha,\iota,\Upsilon}, h^{q_{j'}}_{\alpha,\iota,\Upsilon})
\mbox{ while }
(N^r_{\alpha,\iota,2\Upsilon+1}, h^r_{\alpha,\iota,2\Upsilon+1})=
(N^{q_j}_{\alpha,\iota,\Upsilon}, h^{q_j}_{\alpha,\iota,\Upsilon}).
\]
We have now completed the proof of Subcase B of Case 2 of the Claim.
\smallskip
\underline{Subcase A}. $\lambda=\aleph_0$. We have to prove that $R$ satisfies
ccc. Let
functions $g^\ast_0, g^\ast_1, g^\ast_2$ and $g^\ast_3$ be as in the proof of
Subcase B, and let the function $F$ and the club $C$ be given as in that proof.

Suppose that we are given a sequence $\langle q_s:\,s<\omega_1\rangle$ of
conditions in $R$. Let
$E\subseteq C$ be a club of $\omega_1$ such that 
\[
j\in S^{\omega_1}_{\aleph_0}\,\,\&\,\,j'<j\implies(\forall i<i(q,j'))
(w[M^{q_{j'}}_{\eta(q_{j'},i)}]\subseteq j).
\]
Let $g$ be a function that to an ordinal $j\in S^{\omega_1}_{\aleph_0}$
assigns
\[
(g^\ast_3(q_j), \langle w[M^{q_j}_{\eta(q_j,i)}]\cap j:\,i<i(q_j)\rangle)
\]
and let
\[
f=(F\circ g)\rest (E\cap S^{\omega_1}_{\aleph_0})\cup 0_{\omega_1\setminus
(E\cap S^{\omega_1}_{\aleph_0})}.
\]
Exactly as in Subcase B, it follows that whenever
\[
j'<j\in S^{\omega_1}_{\aleph_0}\cap E\,\,\&\,\,\,f(j')=f(j)
\]
then letting
\[
u=u^{q_{j'}}\cup u^{q_j}\cup\{\eta\cap\nu:\,\eta\in u^{q_j},\nu\in u^{q_j}\}
\]
and for $\eta\in u$
\[
M_\eta=M^{q_l}_\eta\rest\lambda\cdot\llg(\eta)
\]
for any $l\in \{j,j'\}$ for which $\eta\in u^{q_l}$, we obtain a well defined
sequence $\langle M_\eta:\,\eta\in u\rangle$ of elements of $K$ with the
property that for any $l\in\{j,j'\}$, $\eta\in u^{q_l}$ we have
$M^{q_l}_\eta\le M_\eta$. Let $S\subseteq S^{\omega_1}_{\aleph_0}\cap E$
be stationary such that $f(j)$ is fixed for $j\in S$. We apply the
$\Delta$-system lemma to $\{b^{q_j}:\,j\in S\}$ and obtain $A\in
[S]^{\aleph_1}$ and
$b^\ast\in [\mu]^{<\aleph_0}$ such that for every $j\ne j'\in A$ we have
$b^{q_j'}\cap b^{q_j}=b^\ast$. If $b^\ast=\emptyset$, then for every $j,j'\in A$,
the condition 
\[
\langle u,\langle M_\eta:\,\eta\in u\rangle, b^{q_{j'}}\cup
b^{q_j},
\langle c^{q_l}_\alpha:\,\alpha\in b^{q_l}\rangle,
\langle (N^{q_l}_{\alpha,\iota},h^{q_l}_{\alpha,\iota}):\,\alpha\in b^{q_l},
\iota\in c^{q_l}_\alpha\rangle\rangle
\]
is a common upper bound of $q_{j'}$ and $q_j$ where $u, \langle M_\eta:\,\eta\in
u\rangle$ are defined above.

Suppose that $\card{b^\ast}=n^\ast>0$. Using the $\Delta$-system lemma $n^\ast$ times
if necessary, we can find $B\in [A]^{\aleph_1}$ and for $\alpha\in b^\ast$ a
set $c^\ast_\alpha\in [\omega_1]^{<\aleph_0}$ such that
$\alpha\in b^\ast\,\,\&\,\,j'<j\in B\implies$
\begin{description}
\item{(i) $c^{q_{j'}}_\alpha\cap c^{q_j}_\alpha= c^\ast_\alpha$,
\item{(ii)} $\min(c^{q_{j'}}_\alpha\setminus
c^\ast_\alpha)>\max\{\lambda\iota:\,\iota\in c^\ast_\alpha\}$,
\item{(iii)} $\min(c^{q_j}_\alpha\setminus
c^\ast_\alpha)>\max\{\lambda\iota:\,\iota\in c^{q_{j'}}_\alpha\}$,
\item{(iv)} $\iota\in c^\ast_\alpha\implies
(N^{q_{j'}}_{\alpha,\iota},h^{q_{j'}}_{\alpha,\iota} )=
(N^{q_j}_{\alpha,\iota},h^{q_j}_{\alpha,\iota} )$
}
\end{description} 
\underline{and} for $k<n^\ast_\alpha\deq\card{c^\ast_\alpha}$
letting $\iota'_k,\iota_k$ be the $k$-th element
of $c^{q_{j'}}_\alpha, c^{q_j}_\alpha$ respectively, we have that
$N^{q_{j'}}_{\alpha,\iota'_k}$ and $N^{q_j}_{\alpha,\iota_k}$ are isomorphic.
Let $j'<j\in B$ and let $\alpha\in b^\ast$. Let $N_j'=\bigcup_{\iota\in
c^{q_{j'}}_\alpha}N^{q_{j'}}_{\alpha,\iota}$ and
$N_j=\bigcup_{\iota\in
c^{q_j}_\alpha}N^{q_j}_{\alpha,\iota}$, while $h_j'=\bigcup_{\iota\in
c^{q_{j'}}_\alpha}h^{q_{j'}}_{\alpha,\iota}$ and
$h_j=\bigcup_{\iota\in
c^{q_j}_\alpha} h^{q_j}_{\alpha,\iota}$. Then $N_j$ and $N_j'$ are isomorphic
and $h_j$ and $h_j$ agree on their intersection, while there are
$\delta_0<\delta_1<\delta_2$ divisible by $\lambda$ such that 
\[
N_{j'}\rest\delta_0=N_j\rest\delta_1=N_j\cap N_j'
\]
with $\card{N_j'}\subseteq\delta_1$ and $\card{N_j}\subseteq\delta_2$
and $\eta=f^\ast_\alpha\rest\delta_2\in u^{q_j}$. We also have that $N_j',
N_j\in \Gamma$. Then $h_j$ is a lawful embedding of $N_j$ into $M_\eta$ and
$h_j'$ is a lawful embedding of $N_j'$ into $M_\eta$, by the choice of $S$. 
Similarly to the proof of Claim \ref{facts}(4), we can see that $q_j'$ and $q_j$
are compatible, by finding $N\in\Gamma$ with $N\ge N_1,N_2$, extending
$q_j$, $q_{j'}$ to enlarge $M_\eta$ and then taking an upper bound of
the extensions.

\underline{Case 3}. Suppose that
\[
\langle (\langle q^\xi_i:\,i<\lambda^+\rangle, f_\xi),
\langle p^\xi_i:\,i<\lambda^+\rangle:\,\xi<\zeta\rangle
\]
have been played so
far. By Claim \ref{A} we can have $q^\zeta_i$ be the lub of
$\langle q^\xi_i:\,\xi<\zeta\rangle$. Let $\langle
I_\gamma:\,\gamma<\lambda\rangle$ list the isomorphism types of elements of
$K$.
Let $F$ be a bijection
\[
F:\,\lambda\times K\times\{h:\,h\mbox{ a lawful function with }
\Dom(h)\in [\lambda^+]^{<\lambda}\}\into\lambda^+.
\]
Let $C$ be a club of $\lambda^+$ such that for
every $j\in S^{\lambda^+}_\lambda\cap C$ we have
\[
\gamma<\lambda,u\in [j]^{<\lambda}\,\,\&\,\,\Dom(h)\in [j]^{<\lambda}
\implies F((\gamma,u,h))<j.
\]
Let $C_\zeta=C$ and define $g_\zeta$ which to an ordinal $j\in
S^{\lambda^+}_\lambda$ assigns 
\[
({\rm type}(N^{q^\zeta_j}), N^{q^\zeta_j}\cap
j, h^{q^\zeta_j}\rest (N^{q^\zeta_j}\cap
j)).
\]
Then let $f_\zeta=(F\circ g_\zeta)\rest(C_\zeta\cap
S^{\lambda^+}_\lambda)\cup 0_{\lambda^+\setminus (C_\zeta\cap
S^{\lambda^+}_\lambda)}$. Let $E\subseteq C$ be a club of $\lambda^+$
such that 
\[
j\in E\cap S^{\lambda^+}_\lambda\,\,\&\,\,j'<j\implies
(\forall\zeta<j')(\card{N^{q^\zeta_{j'}}}\subseteq j).
\]
Let $(N',h')$ be the lub of $\{(N^{q^\zeta_{j'}},
h^{q^\zeta_{j'}}):\,\zeta<\varepsilon\}$ and $(N,h)$ the lub of
$\{(N^{q^\zeta_j},
h^{q^\zeta_j}):\,\zeta<\varepsilon\}$. We shall show that $(N,h)$ and $(N,h')$
are compatible. As $N,N'\in \Gamma_\alpha$, clearly they are compatible as
elements of $K$. We need to show that $h$ and $h'$ agree on $N\cap N'$.

Suppose not and let $\zeta<\varepsilon$ be the least such that $h_\zeta$ and
$h'_\zeta$ disagree on $N^{q^\zeta_{j'}}\cap N^{q^\zeta_j}$-such a $\zeta$
exists by the definition of the lub in the forcing. By the choice of $E$ we
have $\card{N^{q^\zeta_{j'}}}\subseteq j$ and by the choice of $f_\zeta$ we have
$(N^{q^\zeta_{j'}}\cap j', h^{q^\zeta_{j'}}\rest (N^{q^\zeta_{j'}}\cap j'))=
(N^{q^\zeta_j}\cap j, h^{q^\zeta_j}\rest (N^{q^\zeta_j}\cap j))$.
Hence $N^{q^\zeta_{j'}}\cap N^{q^\zeta_j}\subseteq j'$ and $h^{q^\zeta_{j'}}$
and $h^{q^\zeta_j}$ agree on this intersection, a contradiction. We can also
see that $N'$ and $N$ are isomorphic.

Now let $G$ be as in the Definition \ref{candidate} Case 3 and let $p'\in G$
witness that $(N',h')\in R$, while $p\in G$ witnesses that $(N,h)\in R$. Let
$p^+\in G$ be a common upper bound of $p$ and $p'$.
By Claim \ref{facts}(4)
it follows that 
\[
\DD=\{p^{++}\ge p^+:\,(\exists \iota\in c^{p^{++}}_\alpha)(\exists\Upsilon\in 
d^{p^{++}}_{\alpha,\iota}) \,(N,N'\le N^{p^{++}}_{\alpha,\iota,\Upsilon}
\,\,\,\&\,\,h\cup h'\subseteq h^{p^{++}}_{\alpha,\iota,\Upsilon})\}
\]
is dense in the forcing $R'$ giving rise to $G$, which suffices.

Now suppose that
we are in \underline{Case 4} and that \[
\left\langle(\langle q^\xi_i:\,i<\lambda^+\rangle, f_\xi),
\langle p^\xi_i:\,i<\lambda^+\rangle:
\,\xi<\zeta\right\rangle
\]
have been played so far in the game $\ast^\varepsilon_\lambda[R]$.
This is where we get to use the workability of $K$.
As before, we let player I choose $q^\zeta_i$ as the unique least upper bound of
$\langle p^\xi_i:\,\xi<\zeta\rangle$. Let $q^\zeta_i=\langle
M^\zeta_i,N^\zeta_i\rangle$.
Using $\lambda=\lambda^{<\lambda}$,
we can find a
regressive function $f_\zeta$
such that if $i<j$ in $S^{\lambda^+}_\lambda$ are such that
$f_\zeta(i)=f_\zeta(j)$, then 
\begin{description}
\item{(a)} $N^\zeta_i\rest i= N^\zeta_j\rest j$,
\item{(b)} There is a $K_{\ap}$-isomorphism $h^\zeta_{i,j}$
from $N^\zeta_i$ onto $N^\zeta_j$ mapping $M^\zeta_i$ onto
$M^\zeta_j$, and such that $h^\zeta_{i,j}\rest (\card{N_i^\zeta}\cap i)$
is the identity.
\end{description}
At the end,
let $C\subseteq\lambda^+$
be a club such that for every $\zeta<\varepsilon$  
\[
i<j\,\,\&\,\,j\in C\implies \card{N^\zeta_i}\subseteq j.
\]
Suppose now that
$i<j\in C\cap S^{\lambda^+}_\lambda$ are such that $f_\zeta(i)=
f_\zeta(j)$ for all $\zeta<\varepsilon$.

For $l\in \{i,j\}$ let $M_l\deq\bigcup_{\zeta<\varepsilon} M^\zeta_l$
and $N_l\deq\bigcup_{\zeta<\varepsilon} N^\zeta_l$. Notice that
$\langle M_l,N_l\rangle\in R$ and that for every $\zeta<\varepsilon$ we
have $\langle M^\zeta_l,N^\zeta_l\rangle \le \langle M_l, N_l\rangle$. 
Also observe that $\card{N_i}\subseteq j$ and that $N_i\rest i= N_j\rest j$.
Let $h=h_{i,j}\deq\bigcup_{\zeta<\varepsilon}h^\zeta_{i,j}$. Then $h$
is a $K_{\ap}$-isomorphism from $N_i$ onto $N_j$ mapping $M_i$ onto
$M_j$, and such that $h\rest(\card{N_i}\cap i)$ is the identity.
By the definition of workability, we can find
$\langle M,N\rangle \in R$ which is stronger than both
$\langle M_i, N_i\rangle$ and
$\langle M_j, N_j\rangle$.

$\eop_{\ref{astepsilon}}$
\end{Proof of the Claim}

\begin{Claim}\label{whatget} Suppose that $R$ is given by Case 4 of Definition
\ref{candidate}.
Then, keeping the notation of Def \ref{candidate}, in 
$V^{P_\beta\ast\name{Q}^\beta_j\ast\name{R}^\beta_j}$ we have
\[
\Gamma
=\{N':\,(\exists N\ge N')(\exists M)[\langle M,N\rangle \in
G_{R}] \}\in
K_{\md}[K_{\ap}]\},
\]
where $G_R$ is $R$-generic over $V^{P_\beta\ast\name{Q}^\beta_j}$.
\end{Claim}

\begin{Proof of the Claim} We verify that $\Gamma$
satisfies the required properties (i)-(v) from Definition \ref{md}.
As (v) is obvious, we check (i)-(iv). By Claim 
\ref{A}(2)
\[
R \mbox{ is a } (<\lambda)\mbox{-complete
forcing},
\]
hence $\Gamma$ is $(<\lambda)$-closed, and so it satisfies (i). Property (ii)
follows by genericity.

Given $\beta<\lambda^+$ such that the requirement of (iii)
of Definition \ref{md}(2) holds for $\Gamma^-$ ($\Gamma^-$ comes from the
definition of $R$ by Case 4),
arguing in $V^{P_\beta\ast\name{Q}^\beta_j}$ we shall show that
\[
\II\deq\{\langle M,N \rangle \in R:\,(\exists\gamma\in\card{N})[\beta+\lambda=
\gamma+\lambda]\}
\]
is dense in $R$. So let $\langle M,N\rangle \in R$ be given. Let $M'\in \Gamma^-$
be $\ge M$ and such that for some $\gamma$ with $\gamma+\lambda=\beta+\lambda$ we have
$\gamma\in\card{M'}$.
Since $N\rest{\rm Ev}=M$, we can apply amalgamation to $N,M,M'$ to find
$N'\ge_{K_{\rm ap}} N$ with $M'\le N'$.
By Remark \ref{concern} (3), we can assume that
$N'\rest \Ev=M'$. Hence $\langle M',N'\rangle 
\in R\cap\II$ is as required, showing (iii).

To show (iv), suppose $N\in\Gamma$ and $N\le_{K_{\rm ap}} N'$ after we have
forced by $R$.
As the forcing
with $R$ is $(<\lambda)$-closed, we have $N'\in V^{P^\beta\ast\name{Q}^\beta_j}
$ and $N\le_{K_{\rm ap}} N'$ holds. Let $M$ be such that $\langle M,N\rangle \in G_{R}$.
Now observe that by amalgamation and Remark \ref{concern}(3), the set
\[
\{\langle M,N''\rangle:\,(\exists \mbox{ lawful }h)[h:\,N'\into N''\mbox{ embedding over } M]\}
\]
is dense in $R$ above $\langle M,N\rangle$.
$\eop_{\ref{whatget}}$
\end{Proof of the Claim}

This finishes the inductive proof.

\begin{Claim}\label{bookkeeping} It is possible to define the iteration
$\bar{P}$ so that in $V^{\bar{P}}$ we have
\begin{description}
\item{(1)} If $\lambda>\aleph_0$ then for every $(<\lambda)$-complete forcing notion $Q$
which
satisfies $\ast^\varepsilon_\lambda$ and has the set of
elements some ordinal $<\kappa$ and for every $\beta<\lambda^{++}$ large
enough
we have $\name{R}^\beta_j=Q$ for some $j<\mu$. If $\lambda=\aleph_0$, the analogous
statement holds with ccc forcing in place of $(<\lambda)$-complete $\ast^\varepsilon_\lambda$
forcing,

\item{(2)} For every workable
strong $\lambda$-approximation family $K$
and a family $\bar{\Gamma}=\{\Gamma_\alpha:\,\alpha<\mu\}$
of elements of $K_{\rm md}$, and for every $\beta<\lambda^{++}$ large
enough,
there is $j<\mu$ such that $R^\beta_j$ is given by Case 2 of Definition
\ref{candidate}  using $K,\Gamma^-$ as parameters.

\item{(3)} If $\lambda>\aleph_0$, then for every
$K,\bar{\Gamma},\beta,j$ as in (2), for every $\alpha<\mu$,
there is $\beta'>\beta$ such that $R^{\beta'}_j$ is defined by Case 3 of
Definition \ref{candidate} using $\Gamma_{\alpha}$ and $\beta$ as parameters.

\item{(4)} For every workable strong $\lambda$-approximation family $K$
and $\Gamma^-\in K_{\md}^-[K]$
such that $\bigcup\{\card{M}:\,M\in\Gamma\}\subseteq \Ev$,
for every $j$ large enough
there is $\beta<\lambda^{++}$ such that $R^\beta_j$ is defined by Case 4 of
Definition \ref{candidate} using $\Gamma^-$ as a parameter.

\end{description}
\end{Claim}

\begin{Proof of the Claim} We use the standard bookkeeping. As the forcing is
$(<\lambda)$-closed,
any workable strong $\lambda$-approximation family $K\in V^{\bar{P}}$
appears at some stage and does not gain any new members later. Also notice that
being in $K_{\rm md}^-$ and $K_{\rm md}$ is absolute between $V^{\bar{P}}$ and
$V^{P^\beta_j\ast\name{Q}^\beta_j}$ containing $\Gamma$ for $\Gamma\subseteq K$.
$\eop_{\ref{bookkeeping}}$
\end{Proof of the Claim}
$\eop_{\ref{universal}}$

This finishes the proof of the Theorem.
\end{Proof of the Theorem}

\begin{Remark} Applying the usual proof of the consistency of $MA+\neg CH$ if
we assume in Theorem \ref{universal} that $V$ satisfies
\[
\theta<\kappa\implies\theta^{<\lambda}<\kappa,
\]
we can drop the assumptions $\card{Q}<\kappa$ from (d) in the
conclusion
of Theorem \ref{universal}.
\end{Remark}

\begin{Proof of the Conclusion} 
Let $V^\ast\deq V^{\bar{P}}$, where $\bar{P}$ is an iteration satisfying the
requirements listed in Claim \ref{bookkeeping}.

{\noindent (1)}
Given an abstract elementary class $\KK$ in $V^\ast$
such that there is a workable strong $\lambda$-approximation
family
$K_{\ap}$ approximating $\KK$ and such that
${\rm LS}(\KK)\le\lambda$ and suppose that $M\in \KK_{\lambda^+}$. Let
$\langle \bar{\Delta}_\beta:\,\beta<\lambda^{++}\rangle$ be as in (e) of the
conclusion of Theorem \ref{universal}, for our $K_{\ap}$. Let $M_\beta^\ast$
be as
in Claim \ref{tree}, with $M_\eta\deq\bigcup\Delta^\beta_\eta$ for
$\eta\in\TT$. (Note that $\eta\initialeq\nu$ does imply that $M_\eta\le_{\KK_{\ap}}M_\nu$).
We claim that $M$ embeds into $M^\ast_\beta$ for some $\beta$.

By the definition of approximation, there
is $\Gamma^-$ which
is an element of $K^-_{\md}[K_{\ap}]$, such that $M\le_{\KK} M_{\Gamma^-}$
and $N\in\Gamma^-\implies\card{N}\subseteq{\rm Ev}.$
By Theorem \ref{universal}(f), there is $\Gamma\in K_{\md}[K_{\ap}]$
such that $\Gamma^-\subseteq \Gamma$, and
hence by Observation \ref{union}(2), we have $M_{\Gamma^-}
\le_{\KK} M_\Gamma$. Let $\beta<\lambda^{++}$ be such that $\Gamma, K_{\rm ap}
\in V_\beta$, which is easily seen to exist.
By (e) in the conclusion of Theorem
\ref{universal} and its proof, there is $\alpha<\mu$
such that $M_\Gamma$ is
isomorphically embeddable into
$M_{\bigcup_{i<\lambda^+}\Delta^\beta_{f^\ast_\alpha\rest i}}$. 
By Observation \ref{dodatna}, we have
\[
M\le_{\KK}M_{\Gamma^\ast}\le_{\KK}M_\Gamma\le_{\KK}
M_{\bigcup_{i<\lambda^+}\Delta^\beta_{f^\ast_\alpha\rest i}}
\le_{\KK}M^\ast_\beta.
\]

{\noindent (2)} In addition to what we have already observed, we need
to observe that $2^\lambda=\mu$, and this is the case because
$\bar{P}$ adds a Cohen subset to $\lambda$ $\mu$ many times.

{\noindent (3)} Follows from (1) of the Theorem.

{\noindent (4)} This part follows similarly to (1), using the assumptions
on $\KK^+$.
$\eop_{\ref{conclus}}$
\end{Proof of the Conclusion}

\begin{Fact}\label{whichones} Suppose $\lambda=\lambda^{<\lambda}
\ge\aleph_0$.
Each of the following classes $\KK$ is an abstract elementary class
for which there is a workable strong $\lambda$-approximation family
approximating it,
and the L\"owenheim-Skolem number of $\KK$ is $\le\lambda$:

\begin{description}
\item{(1)} The class
of models of $T^\ast_{\rm feq}$, i.e. an indexed family of
independent equivalence relations, with $M\le N$ iff $M$ embeds
into $N$,
\item{(2)} The class $T_{trf}$ of triangle free graphs, with the
same order as in (1),
\item{(3)} The class of models of any simple theory.
\end{description}
\end{Fact}

[Why? (1) and (2) were proved in [Sh 457], and (3) is proved in
\cite{Sh500}.]

\section{Consistency
of the non-existence of universal normed vector spaces.}\label{negativex}

\begin{Definition}\label{Banach}  Suppose that $I$ is a linear order.

{\noindent (1)} 
We define
a vector space $B_I$ over ${\Bbf Q}$
by
\[
B_I\deq\left\{\sum_{i\in I} a_i x_i:\,
a_i\in \Bbf{Q}\,\,\&\,\,\{i:\,a_i\neq 0\}\mbox{ finite }\right\},
\]
where $\{x_i:\,i\in I\}$ is a set of variables that serve as a
basis for $B_I$.
The addition and scalar multiplication is defined in the obvious manner.

{\noindent (2)}
For any
$I$-increasing sequence $\bar{t}\in {}^{\omega>} I$, we define
a functional $f_{\bar{t}}:\,B_I\to\Bbf{R}$ by letting
\[
f_{\bar{t}}(\sum_{i\in I} a_i x_i)\deq \sum_{l<\lg(\bar{t})}
{{1}\over{\ln(l+2)}}a_{\bar{t}(l)}.
\]
Let
\[
F\deq\{f_{\bar{t}}:\,\bar{t}\in\bigcup_{n<\omega}{}^n I\mbox{ is }
I\mbox{-increasing}\}.
\]
For $x\in B_I$ we define
$\norm{x}=\norm{x}_F\deq\sup\{\card{f(x)}:\,f\in F\}$.

\end{Definition}

\begin{Note}\label{jesuli}
(1) Functionals $f_{\bar{t}}$ defined as above are linear.

{\noindent (2)} For every $x\in B_I$, there are only finitely many possible values of
$f_{\bar{t}}(x)$. (Hence, $\norm{x}
=\Max\{\card{f(x)}:\,f\in F\}$).
\end{Note}

\begin{Claim}\label{defined} Suppose that $I$ and $B_I$
are as in Definition \ref{Banach} and $I$
is infinite. Then
$B_I$ is a normed vector space over $\Bbf{Q}$ with $\card{B_I}
=\card{I}$.
\end{Claim}

\begin{Proof of the Claim} We prove that
$\norm{-}$ is a norm on $B_I$.
The triangular inequality is easily verified. We need to check that
for all $x\in B_I$ we have $0\le \norm{x}<\infty$ and
$\norm{x}=0\iff x=0$.
The second statement is obvious, by considering
sequences $\bar{t}$ whose length is 1, and
the first follows from Note \ref{jesuli}(2).
$\eop_{\ref{defined}}$
\end{Proof of the Claim}

\begin{Theorem}\label{negative} Suppose that $\aleph_0\le
\lambda=\lambda^{<\lambda}
<\lambda^+<\mu=\cf(\mu)=\mu^{\lambda^+}$.

Then for some $(<\lambda)$-complete
and $\lambda^+$-cc forcing notion $P$
of cardinality $\mu$, we have that $P$
forces
\begin{description}
\item{}
$``2^\lambda=\mu$ and for every normed vector
space
$\name{A}$ over $\Bbf{Q}$ of cardinality $\card{\name{A}}<\mu$, there is a
normed vector space $\name{B}$ over $\Bbf{Q}$ of dimension $\lambda^+$
(so cardinality $\lambda^+$) such that there is no vector
space embedding $h:\,\name{B}\to \name{A}$ with the property that for some
$\name{c}\in \Bbf{R}^+$ for all $x\in \name{B}$
\[
1/\name{c}<{{\norm{h(x)}_{\name{A}}}\over{\norm{x}_{\name{B}}}}<\name{c}."
\quad\quad\quad\quad\quad(\ast)
\]
\end{description}
\end{Theorem}

\begin{Proof} We deal with the situation $\lambda>\aleph_0$,
and the proof for $\lambda=\aleph_0$ is similar but easier.

\begin{Definition} (1) We define an iteration
\[
\langle P_\alpha,\name{Q}_\beta:\,\alpha\le\mu,\beta<\mu\rangle
\]
with $(<\lambda)$-supports such that for all $\beta<\mu$ we
have that $\name{Q}_\beta$ is a $P_\beta$-name defined by
\[
\name{Q}_\beta\deq\{(\name{w},\name{\le}_w):\,\name{w}\in
[\lambda^+]^{<\lambda}\,\,\&\,\,
\name{\le}_w\mbox{ is a linear order on }\name{w}\},
\]
ordered by letting $(\name{w},\name{\le}_w)\le (\name{z},\name{\le}_z)$
iff $\name{w}\subseteq \name{z}$ and $\name{\le}_w=
\name{\le}_z\rest(\name{w}\times\name{w})$.

{\noindent (2)} Let $P\deq P_{\mu}$.
\end{Definition}

\begin{Claim}\label{cccc} (1) For every $\alpha<\mu$ we have
\[
\forces_{P_\alpha}``\name{Q}_\alpha
\mbox{ is }(<\lambda)\mbox{-complete and }
\mbox{satisfies }\ast^\omega_\lambda".
\]
(2) $P$ is $\lambda^+$-cc, $(<\lambda)$-complete and
$\forces_P``2^\lambda=\mu".$

{\noindent (3)} For $\alpha<\mu$ 
\[
\name{I}_\alpha\deq(\lambda^+,\bigcup\{\name{\le}_w:\,
(\name{w},\name{\le}_w)\in \name{G}_{\name{Q}_\alpha}\})
\]
is a $P_{\alpha+1}$-name which is forced
to be a linear order
on $\lambda^+$.
\end{Claim}

\begin{Proof of the Claim} (1) The first statement is obvious, we shall prove the
second one.
The proof is by induction on $\alpha$.
Given $\alpha<\mu$, by the induction hypothesis we have
in $V^{P_\alpha}$ that $\lambda^{<\lambda}=\lambda$. We work in
$V^{P_\alpha}$,
and describe the winning strategy of player I in the game
$\ast^\omega_\lambda[Q_\alpha]$. As $\lambda^{<\lambda}=\lambda$, we can fix
a bijection $F$ which to every triple $(w,\le,\gamma)$, where
$w\in [\lambda^+]^{<\lambda}$ and $\le$ is a linear order on $w$, and
$\gamma<\lambda^+$, assigns an element of $\lambda^+$.
We can find a club $E$ of $\lambda^+$ such that for every $j\in S^{\lambda^+}_\lambda
\cap E$ and every relevant triple $(w,\le,\gamma)$,
\[
w\in [j]^{<\lambda}\,\,\&\,\,\gamma<j\implies F((w,\le,\gamma))<j.
\]

Suppose that $n<\omega$ and
\[
\left\langle\langle q^k_i:\,i<\lambda^+\rangle, f_k,
\langle p^k_i:\,i<\lambda^+\rangle:\,k\le n\right\rangle
\]
have been played so far, and we shall describe how to choose $q^{n+1}_i$
and $f_{n+1}$. We let $q^{n+1}_i\deq p_i^n$, for $i<\lambda^+$.
Suppose that $p_i^n
=(w_i^n,\le_i^n)\in Q_\alpha$ for $i<\lambda^+$.
For $j<\lambda^+$, let $\gamma(j,n)\deq\sup(w^n_j\cap j)$.
Note that for $j\in S^{\lambda^+}_\lambda$ we have $\gamma(j,n)<j$. Let $C_{n+1}
\deq E$. Define $g_{n+1}$ which to an ordinal $j\in S^{\lambda^+}_\lambda$
assigns
\[
(w^n_j\cap j, \le^n_j\rest (w^n_j\cap j), \gamma(j,n)).
\]
Then let 
\[
f_{n+1}\deq (F\circ g_{n+1})\rest (C_{n+1}\cap S^{\lambda^+}_\lambda)
\cup 0_{\lambda^+\setminus (C_{n+1}\cap S^{\lambda^+}_\lambda)}.
\]
Hence $f_{n+1}$ is regressive on $C_{n+1}\setminus \{0\}$.

At the end of the game, for $i<\lambda^+$ let
$w^i\deq\bigcup_{n<\omega} w^n_i$
and $\le^i\deq\bigcup_{n<\omega}\le^n_i$. Let $C\subseteq E$ 
be a club such that $i<j\in C\implies w^i\subseteq j$. Suppose
that $i<j\in C\cap S^{\lambda^+}_\lambda$ are such that for all $n<\omega$
we have $f_n(i)=f_{n}(j)$. We shall define a condition $p$ such that $p=(z,\le_z)$
and $z=w^i\cup w^j$, by amalgamating linear orders. 
For $x,y\in z$ we let $n=n(x,y)$ be the minimal $n$
such that $x,y\in w^n_i\cup w^n_j$,
and let
$x\le_z y$ iff
\begin{description}
\item{(i)} $x,y\in w^l$ and $x\le^ly$ for some $l\in\{i,j\}$, or
\item{(ii)} $x\in w^n_i\setminus w^n_j$ and $y\in w^n_j\setminus w^n_i$
and for some $z\in w^n_j\cap w^n_i$ we have $x\le^i z$ and $z\le ^j y$,
\item{(iii)} $y\in w^n_i\setminus w^n_j$ and $x\in w^n_j\setminus w^n_i$
and (ii) does not hold.
\end{description}
It is easily seen that $p$ is as required.

(2) That $P$ is $\lambda^+$-cc follows from (1) by the fact that
$\ast^\omega_\lambda$ is preserved
under $(<\lambda)$-support iterations. See \cite{Shxx}. That $\forces_P
``2^\lambda=\mu"$ is seen by observing that every $\name{Q}_\alpha$ adds a
subset to $\lambda$.

(3) Obvious.
$\eop_{\ref{cccc}}$
\end{Proof of the Claim}

Suppose that in $V^P$ we have a normed vector space $A$ over $\Bbf{Q}$
with $\card{A}<\mu$, with the universe of $A$ a set of ordinals.
Hence for some $\alpha<\mu$ and
a $P_\alpha$-name $\name{A}$ we have that $A=\name{A}_G$.
Suppose that $h\in V^P$ is a vector space embedding
from $B_{I_\alpha}$ into $A$, satisfying $(\ast)$ above, for some $c\in
\Bbf{Q}$.
Hence for some
$p^\ast\in P/P_\alpha$ we have that
$p^\ast$ forces over $V^{P_\alpha}$ the following
statement:
\[
``\name{h}:\,B_{\name{I}_\alpha}\into A
\mbox{ is a normed vector space embedding satisfying }
(\ast) \mbox{ for }\name{c}."
\]
Without loss of generality, $p^\ast$ decides the
value $c$ of $\name{c}$. Let $0<n^\ast<\omega$ be such that $c<n^\ast$.
Let $\name{x}_i$ for $i<\lambda^+$ be the generators of
$B_{\name{I}_\alpha}$.
For $i<\lambda^+$ we find $p_i\in P/{P_\alpha}$ such that
$p^\ast\le p_i$ and
\[
p_i\forces``\name{h}(\name{x}_i)=y_i"\mbox{ for some }y_i.
\]
Let us now work in $V^{P_\alpha}$. Let $p_i(\alpha)=(w_i,<_i)$,
for $i
<\lambda^+$. Without loss of generality we have $i\in w_i$
for all $i$.

By a $\Delta$-system argument,
noting that $\lambda^{<\lambda}=\lambda$ holds in
$V^{P_\alpha}$, we can find
$Y\in [\lambda^+]^{\lambda^+}$
such that
\begin{description}
\item{(a)} for some $w^\ast$ we have that $w_i\cap w_j=w^\ast$,
for all $i\neq j\in Y$, and $<_i\rest (w^\ast \times w^\ast)$
is constant,
\item{(b)} If $i<j$ are both in $Y$, then
\[
\sup(w_i)<\min(w_j\setminus w^\ast), \mbox{ while }
\sup(w^\ast)<\min(w_i\setminus w^\ast),
\]
\item{(c)} 
$Y\cap w^\ast=\emptyset$,
\item{(d)} For $i<j$ both in $Y$, there is
an isomorphism $h_{i,j}$ mapping
$(w_i,<_i)$ onto $(w_j,<_j)$ such that $h_{i,j}(i)=j$. (Note that $i\in w_i\setminus w^\ast$ and
$j\in w_j\setminus w^\ast$.)
\end{description}

\begin{Observation}\label{convergence} The series
\[
\sum_{l\ge 1}{{1}\over{(l+1)\ln(l+2)}}
\]
diverges, while the sum
\[
\sum_{l=1}^{n}{{1}\over{(n-l+1)\ln(l+2)}}
\]
is uniformly and strictly bounded by 4.
\end{Observation}

\begin{Proof of the Observation} The first statement follows by
comparison
with $\int_{1}^{\infty}{{1}\over{(x+1)\ln(x+2)}}dx$. The second statement
follows from the following estimate:
\[
\begin{array}{l}
\sum_{l=1}^{n}{{1}\over{(n-l+1)\ln(l+2)}}
\le\sum_{l=1}^{[n/2]}{{1}\over{(n-l+1)\ln(l+2)}}
+\sum_{l=[n/2]+1}^{n}{{1}\over{(n-l+1)\ln(l+2)}}\le\\
{{1}\over{[n/2]+1}}\sum_{l=1}^{[n/2]}{{1}\over{\ln(l+2)}}
+\sum_{l=1}^{n-[n/2]}{{1}\over{l\ln(n-l+3)}}\le\\
{{1}\over{[n/2]+1}}\cdot [n/2]\cdot {{1}\over{\ln 3}}
+{{1}\over{\ln([n/2]+2)}}\sum_{l=1}^{n-[n/2]}{{1}\over{l}}\le\\
{{1}\over{\ln 3}}+{{1}\over{\ln([n/2]+2)}}(1+\int_{1}^{n-[n/2]}{{1}
\over{x}}
dx)\le
{{1}\over{\ln 3}}+{{1}\over{\ln([n/2]+2)}}\cdot (1+\ln(n-[n/2]))\\
\le {{1}\over{\ln 3}}+ {{1}\over{\ln 2}}+1 <4.
\end{array}
\]
$\eop_{\ref{convergence}}$
\end{Proof of the Observation}

By Observation 
\ref{convergence}, we can choose $m$ large enough such that
\[
\sum_{l=1}^{m}{{1}\over{(l+1)\ln(l+2)}}\ge 4(n^\ast)^2.
\]
Let us choose
$i_1<\cdots<i_{m}\in Y$.

\begin{Claim}\label{order}
We can find $q'$ and $q''$ in $Q_\alpha$,
both extending all $p_{i_l}(\alpha)$ for $1\le l\le m$, and such that
\[
q'\forces_{Q_\alpha}``\langle i_1,\ldots,i_{m}\rangle
\mbox{ is increasing in }\name{I}_\alpha"
\]
and
\[
q''\forces_{Q_\alpha}``\langle i_1,\ldots,i_{m}\rangle
\mbox{ is decreasing in }\name{I}_\alpha."
\]
\end{Claim}

\begin{Proof of the Claim}
Notice that for no
$1\le l_1<l_2\le m$ and $i\in Y$ do we have that $p_i$
decides the order between $i_{l_1}$ and $i_{l_2}$ in
$I_\alpha$, by the choice of $Y$ (this is elaborated below).
The proof can proceed by induction on $m$. The inductive step
is as in the
proof of $\ast^\omega_\lambda$. The only constraint we could have to letting
$i_{l_1}\le i_{l_2}$ (for $q'$) or $i_{l_2}\le i_{l_1}$ (for $q''$) would be if
some $z\in w^\ast$ would prevent this, but this does not happen.
For example, if we could not let $i_{l_1}\le i_{l_2}$ in $q'$ then this would
mean that $i_{l_1}\ge i_{l_2}$ would have to hold. By the choice of $Y$
and since $i_{l_1}\in w_{i_{l_1}}\setminus w^\ast$ and
similarly for $i_{l_2}$, this could
only be the case if for some $z\in w^\ast$ it would hold that $i_{l_1}\ge_{w_{i_{l_1}}}
z$ while $i_{l_2}\le_{w_{i_{l_2}}}
z$. However, this would contradict item (d) in the choice of $Y$.
$\eop_{\ref{order}}$
\end{Proof of the Claim}

Back in $V^{P_\alpha}$,
let $z\deq \sum_{l=1}^m {{1}\over{l+1}}x_{i_l}$.
Let $a\deq\norm{\sum_{l=1}^m {{1}\over{l+1}}y_{i_l}}_A$.
Hence
\[
q'\forces``\norm{\name{z}}_{B_{\name{I}_\alpha}}\ge 4(n^\ast)^2",
\]
and so $a\ge 4(n^\ast)^2/n^\ast=4 n^\ast$.
On the other hand,
\[
q''\forces``\norm{\name{z}}_{B_{\name{I}_\alpha}} <4",
\]
and hence $a <4{n^\ast}$, a contradiction.
$\eop_{\ref{negative}}$
\end{Proof}
\eject

\eject
\end{document}